\documentclass[10pt,letter]{amsart}

\usepackage{amsmath, amsthm}
\usepackage{amssymb} 
\usepackage{amstext}
\usepackage{amsbsy}
\usepackage{hyperref}

\usepackage{multicol,caption}
\usepackage[dvipsnames]{xcolor}
\usepackage{natbib}
\usepackage{soul}
\usepackage{color}
\usepackage{amscd}

\usepackage{xcolor}
\usepackage{tabularx}
\usepackage{subcaption}
\usepackage{url}
\usepackage{graphicx}
\usepackage{mathtools}

\usepackage{enumerate}
\usepackage{chngpage}

\usepackage{mathtools}

\newtheorem{theorem}{Theorem}[section]
\newtheorem{conjecture}[theorem]{Conjecture}
\newtheorem{question}[theorem]{Question}

\newtheorem{prop}[theorem]{Proposition}
\newtheorem{lemma}[theorem]{Lemma}
\newtheorem{corollary}[theorem]{Corollary}

\newcommand{\fp}{Four$+$}

\newcommand{\Clues}{\texttt{Clues}}

\usepackage{graphicx}

\usepackage{fullpage}
\makeatletter
\newtagform{alphanum}[\thesubsection\alph{equation}\@gobble]{(}{)}
\makeatother

\title{Maximizing Crossword Grids with a Fixed Minimum Clue Length}
\author{Kelly Isham}
\address{Colgate University, Hamilton NY 13346}
\email{kisham@colgate.edu}
\author{Abigail Schnitzer}
\address{Colgate University, Hamilton NY, 13346}
\email{aschnitzer@colgate.edu}
\date{}
\begin{document}
\maketitle

\begin{abstract}
In an American-style crossword puzzle, the grid must be connected, 180-degree symmetric, and have minimum word length 3. Ferland and Pratt construct $n\times n$ crossword grids with many clues satisfying these conditions and prove their construction is maximal when $n\equiv 1$ or 3 $\pmod{4}$. In this paper, we explore American-style crosswords in which the minimum word length is 4. Jacobs computed the maximal $15 \times 15$ grid with minimum word length 4. In this paper, we construct grids with a maximal number of clues for all $n \not \equiv 0 \pmod{5}$.  
\end{abstract}

\section{Introduction}
\label{Intro}

Crosswords have long delighted eager puzzlers and devoted wordsmiths. These puzzles, while seemingly word-oriented, are filled with a myriad of mathematical patterns that continue to be created and discovered. A crossword puzzle's difficulty is directly related to the length and number of clues: puzzles with longer clues are harder since there are fewer additional clues for solvers to rely on.  

In order to delve into crossword counting we must first establish definitions. The $n\times n$ \text{grid} consists of squares colored either black or white. A \textit{clue} is assigned to each maximal linear segment of white squares. For each of these segments the top leftmost white square contains an assigned number which is the clue's label. Therefore, clues with a horizontal answer are listed by their numbers in the Across clues section and  clues with vertical answer are listed by their numbers in the Down clues section. Each letter must be used in both an across and down clue. 

A typical American-Style crossword, as seen for example in the New York Times puzzle section \cite{nytimes}, must follow a few more rules (see e.g., \cite{Cruciverb}). The rules are as follows. 
    \begin{itemize}
  \item \textit{Connectivity}: The grid must have all-over interlock, meaning the black squares may not cut the grid into separate pieces. This ensures that all parts of the puzzle rely on each other, one cannot solve one section on its own. 
  \item \textit{Symmetry}: The grid must look the same when rotated 180\textdegree. This tells us that the pattern of black and white squares is symmetrical. 
  \item \textit{Three+}: The minimum clue length is 3 letters long. 
\end{itemize}

In \cite{ferlandpratt}, Ferland and Pratt consider the maximum number of clues that fit into an $n \times n$ American-style crossword. Ferland also discusses these problems in \cite{ferland17, article}. Ferland poses several open questions about the maximum number of clues in crosswords with relaxed rules in \cite{article}. In his PhD thesis, Bobby Jacobs \cite{jacobs_phd} examines the classic $15 \times 15$ crossword without symmetry and with different word length restrictions. In this paper, we answer the following question, posed by Ferland in \cite{article}.

\begin{question}\cite[Question 3]{article}
    What is the maximum number of clues possible in an American-style crossword puzzle grid so that each clue is at least 4 letters long?
\end{question}

We maintain connectivity and symmetry, but change the minimum clue length to 4.  We call this rule \fp.
In this paper we compute $a(n)$ exactly for all $n \not\equiv 0 \pmod{5}$.

\begin{theorem}\label{THM:main}
Let $a(n)$ be the maximum number of clues in an $n \times n$ American-style crossword with \fp. Let 
\begin{equation*} g(n)= \begin{cases}
    \frac{2}{5}(n^2-3n+20), & n\equiv0 \pmod 5\\
     \frac{2}{5}(n^2-3n+12), & n\equiv1 \pmod 5\\
      \frac{1}{5}(2n^2-7n+14+(n \hspace{-.1in}\pmod{10}), & n\equiv2 \pmod 5\\
   \frac{2}{5}(n^2-3n), & n\equiv3 \pmod 5\\
   \frac{2}{5}(n^2-3n-4)  & n\equiv4 \pmod 5.\\
  \end{cases}
\end{equation*}
    Let $n\geq11$ such that $n \not \equiv 0 \pmod{5}$ and $n\ne 14, 19$. Then $a(n) = g(n)$.
\end{theorem}

Further, we use explicit techniques to find $a(n)$ for $n < 11,$ and $n=14, 19$.

\begin{prop}
  If $4\le n \le 8$, then $a(n) = 2n$. Further, $a(9) = 26$ and $a(10) = 36,  a(14) = 64,$ and $a(19) = 122.$
\end{prop}

Our counts $g(n)$ come from explicit grid constructions. The idea behind the proof is to first show that our constructed $n\times n$ grid, called $G(n)$ has $g(n)$ clues. Then we prove that $a(n) =g(n)$ in all cases $n\not \equiv 0 \pmod{5}$. 

For $n \equiv 0 \pmod{5}$, we construct the grid $G(n)$ for $n = 15, 20, 25$, and 30 and use integer linear programming to show that $g(n) = a(n)$ in these four cases. While our techniques are not strong enough to show that $a(n) = g(n)$ when $n \equiv 0\pmod{5}$ and $n > 30$, we expect this to be true. We present this as a conjecture.

\begin{conjecture}
    For each $n\equiv 0\pmod{5}$ such that $n > 30,$ $a(n)=g(n)$.
\end{conjecture}

\section{Grid Constructions}
\label{grid constructions}

For $4\leq n\leq8$, the all-white grid achieves the maximum number of clues, and $a(n) = 2n$.
When $n=9,10,14$ we find through basic counting methods, detailed later in Section \ref{99proof section}, that $a(9)=26$, $a(10)=36$, and $a(14)=64$. Their maximal grids can be seen in Figures \ref{99mg}, \ref{1010mg}, and \ref{1414mg} respectively. Proofs for $a(10)$ and $a(14)$ are detailed in Appendices \ref{APPA} and \ref{APPB}. 

\begin{figure}[!htb]
\minipage{0.32\textwidth}
\centering
  \includegraphics[width=.65\linewidth]{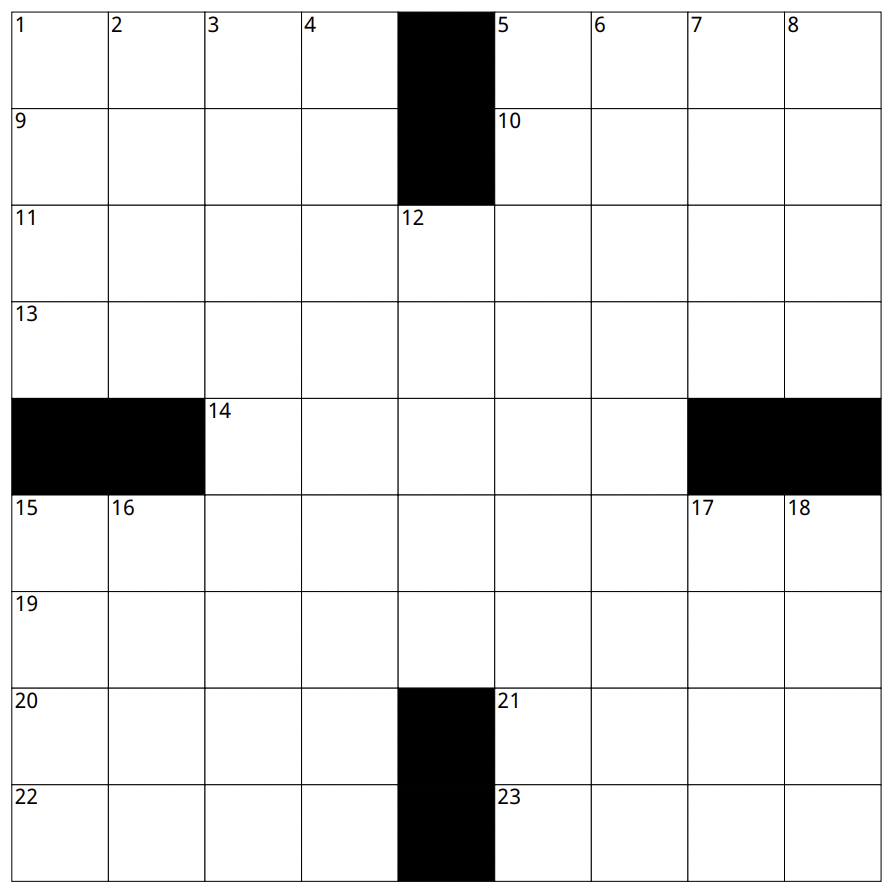}
  \caption{$n=9$ Maximal grid}\label{99mg}
\endminipage\hfill
\minipage{0.32\textwidth}
\centering
  \includegraphics[width=.65\linewidth]{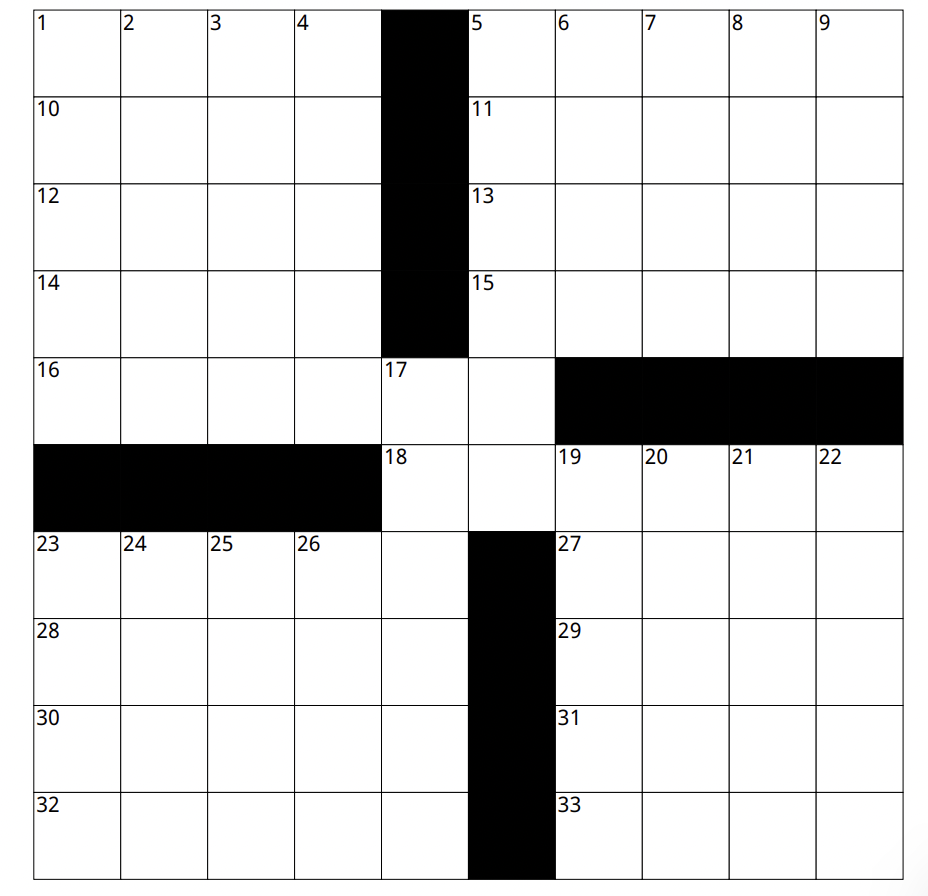}
  \caption{$n=10$ Maximal Grid}\label{1010mg}
\endminipage\hfill
\minipage{0.32\textwidth}%
\centering
  \includegraphics[width=.65\linewidth]{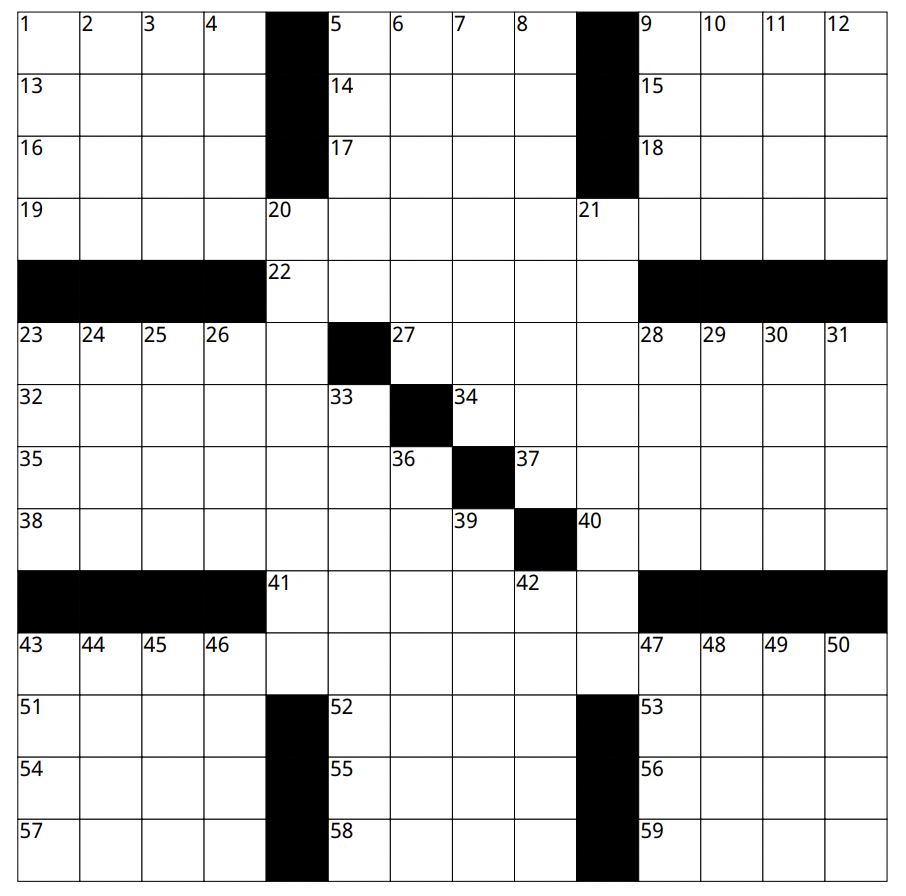}
  \caption{$n=14$ Maximal Grid}\label{1414mg}
\endminipage
\end{figure}

For each $n\geq11$,  we now produce a general construction of a grid $G(n)$. We explicitly compute the number of clues in $G(n)$, denoted $g(n)$. We prove that $g(n)$ is maximal for all $n \ge 11$ and  
$n\not\equiv 0 \pmod{5}$, except for $n=14$ and $n=19$.

\begin{figure}[ht]
    \centering
    \includegraphics[width=0.1\linewidth]{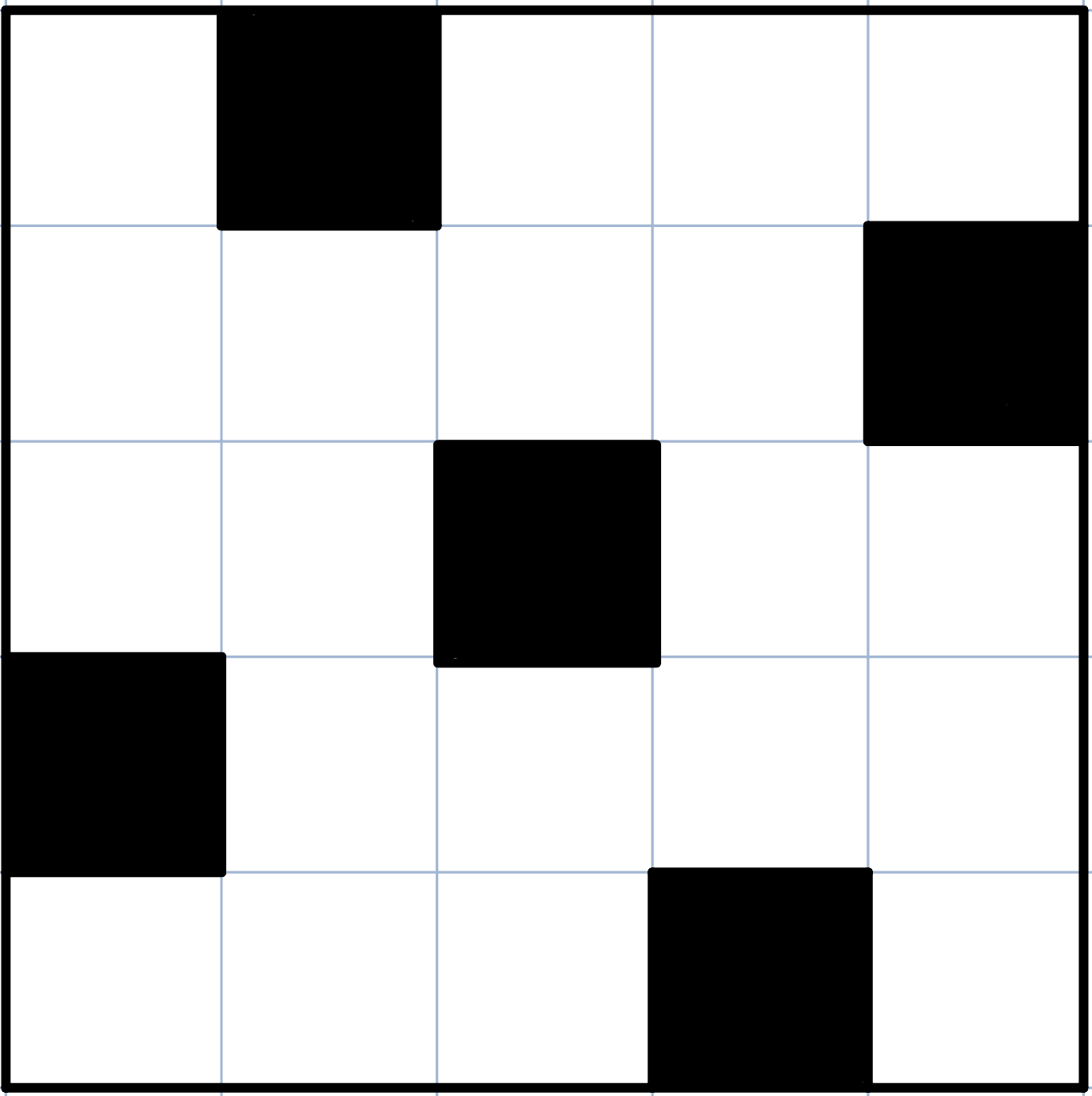}
    \caption{Fundamental Pane}
    \label{fundamental pane}
\end{figure}

The key to the construction is the repeated use of the $5\times 5$ grid shown in Figure \ref{fundamental pane} that we call a \textit{pane}. Note that the fundamental pane has 180 degree rotational symmetry (and in fact, has 90 degree rotational symmetry as well). 

When $n \not \equiv 2 \pmod{5}$, we construct $G(n)$ by first placing our fundamental panes. If $n$ is odd, we place the center square of our fundamental pane in the square of rotation of the $n\times n$ grid. We then tile with fundamental panes until we reach the first/last 3 rows and columns. If $n$ is even, we place the bottom right corner of the fundamental pane at the point of rotation of the $n\times n$ grid and tile as in the odd case. Once we reach the first/last 3 rows and columns, we place the frame. We find the rows and columns for which there is a black square in the 3rd or $(n-2)^{nd}$ row or column and we fill in the above, below, left, or right rows/columns with black squares to maintain Four$+$. 

Note that once this grid is complete, the first/last 4 rows and columns form the frame, which consists of length-4 blocks of black squares. The inner $(n-8) \times (n-8)$ grid consists of tilings of the fundamental pane, following 180-degree symmetry.

If $n \equiv 2 \pmod{5}$, the above construction is not maximal. Instead, we use the fact that there are 4 additional rows in the middle of the $(n-8) \times (n-8)$ inner grid. 

If $n \equiv 2 \pmod{10}$, place into the $\frac{n}{2}$ row, 4 black squares in rows 1 through 4 and maximize this row. By symmetry, this creates a mirror image in the $\frac{n}{2}+1$ row.

If $n \equiv 7\pmod{10}$, place  a 4 row strip into the center, then tile as in the odd case. In this strip, we place black squares at $(\frac{n-1}{2}, 5), (\frac{n-1}{2}, 9) \ldots, (\frac{n-1}{2}, n-8)$ and at $(\frac{n-1}{2}+1, 5), (\frac{n-1}{2}+1, 9) \ldots, (\frac{n-1}{2}+1, n-8)$. Symmetrically, this pattern appears in rows $\frac{n+1}{2}$ and $\frac{n+1}{2}+1$.

\begin{figure}[!htb]
\minipage{0.48\textwidth}
\centering
 \caption*{$18\times18$}
  \includegraphics[width=.50\linewidth]{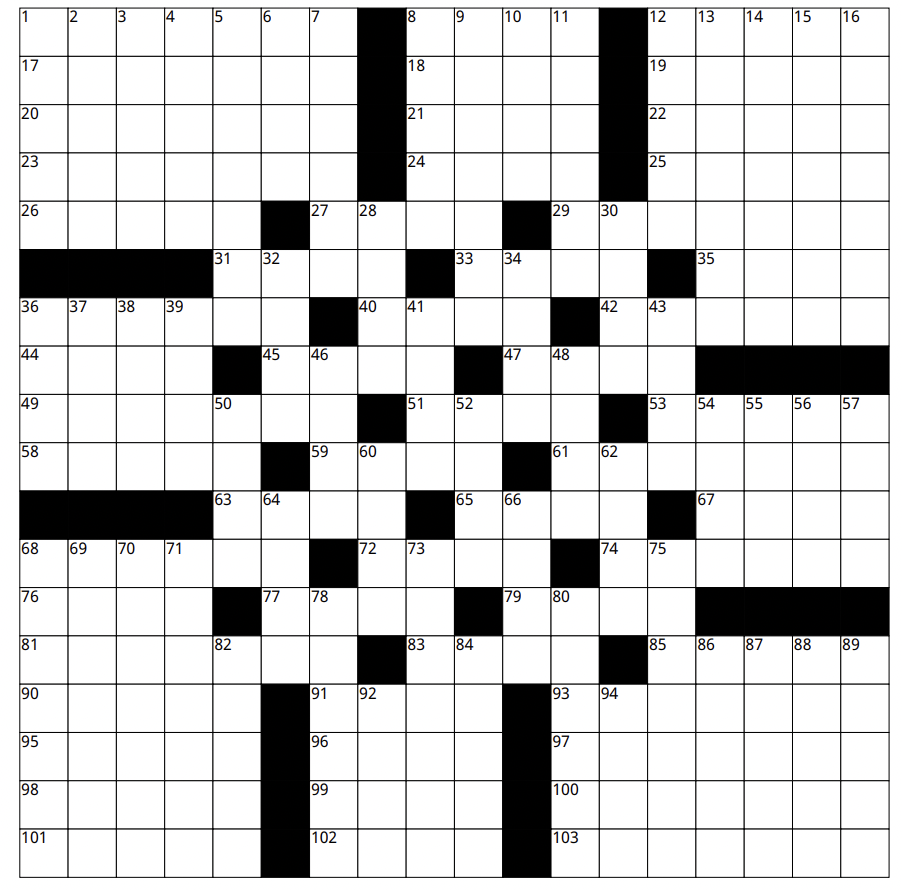}
  \caption{Example of $G(n)$ for $n\equiv3\pmod{5}$}\label{18}
\endminipage\hfill
\minipage{0.48\textwidth}%
\centering
\caption*{$23\times23$}
  \includegraphics[width=.50\linewidth]{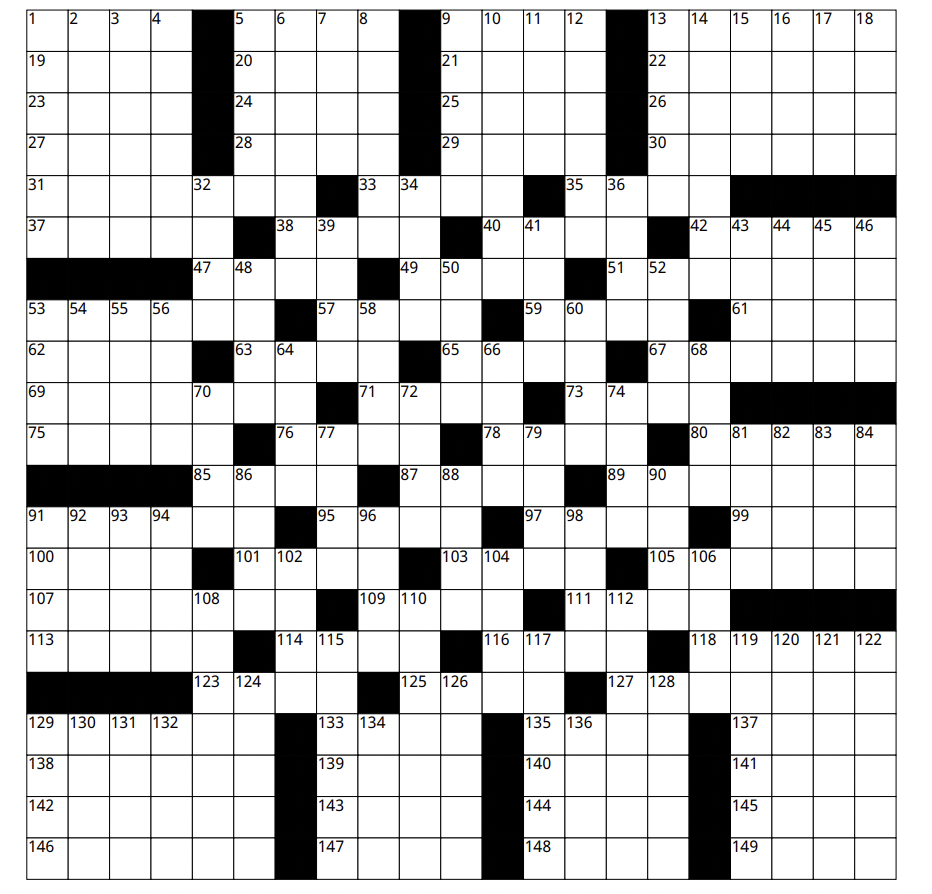}
  \caption{Example of $G(n)$ for $n\equiv3 \pmod{5}$}\label{23}
\endminipage
\end{figure}

\begin{figure}[!htb]
\minipage{0.48\textwidth}
\centering
 \caption*{$22\times22$}
  \includegraphics[width=.50\linewidth]{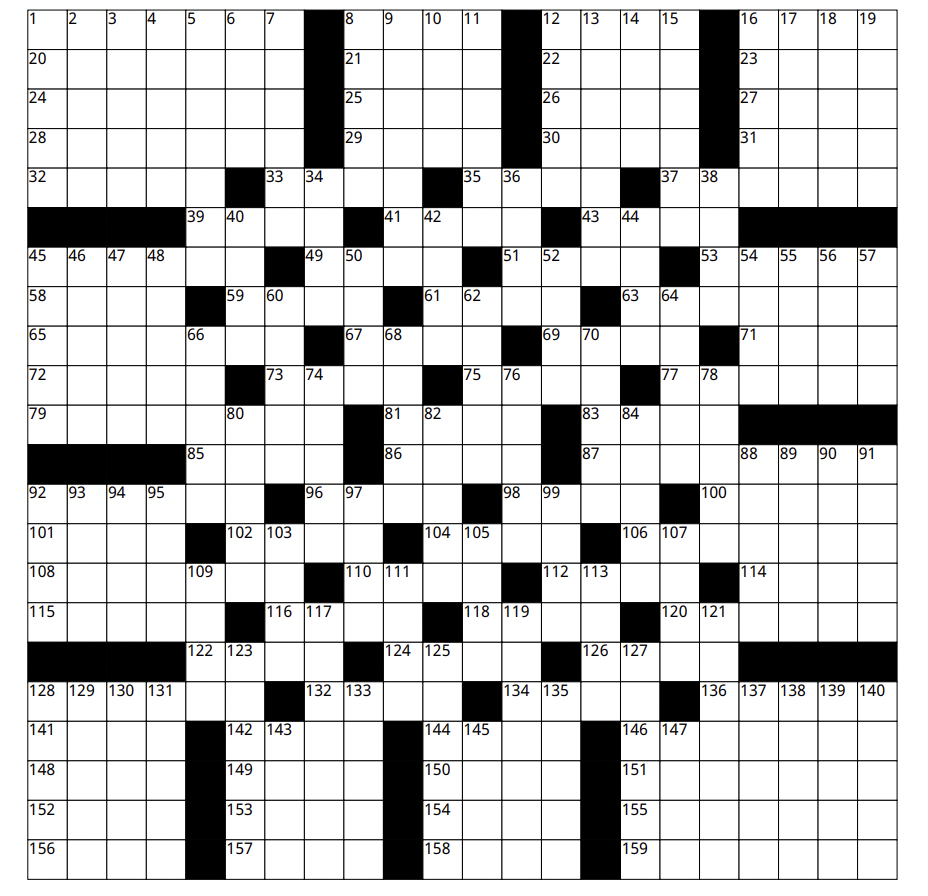}
  \caption{Example of $G(n)$ for $n\equiv2 \pmod{10}$}\label{22}
\endminipage\hfill
\minipage{0.48\textwidth}%
\centering
\caption*{$27\times27$}
  \includegraphics[width=.50\linewidth]{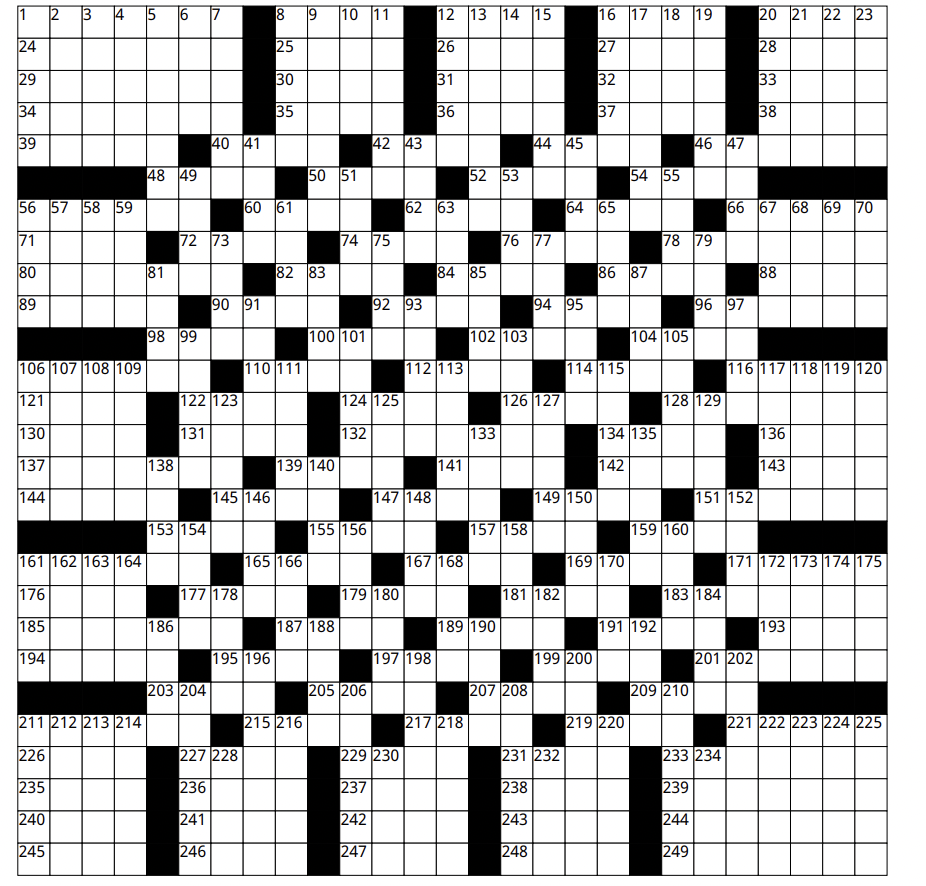}
  \caption{Example of $G(n)$ for $n\equiv7 \pmod{10}$}\label{27}
\endminipage
\end{figure}

\begin{figure}[!htb]
\minipage{0.48\textwidth}
\centering
 \caption*{$24\times24$}
  \includegraphics[width=.50\linewidth]{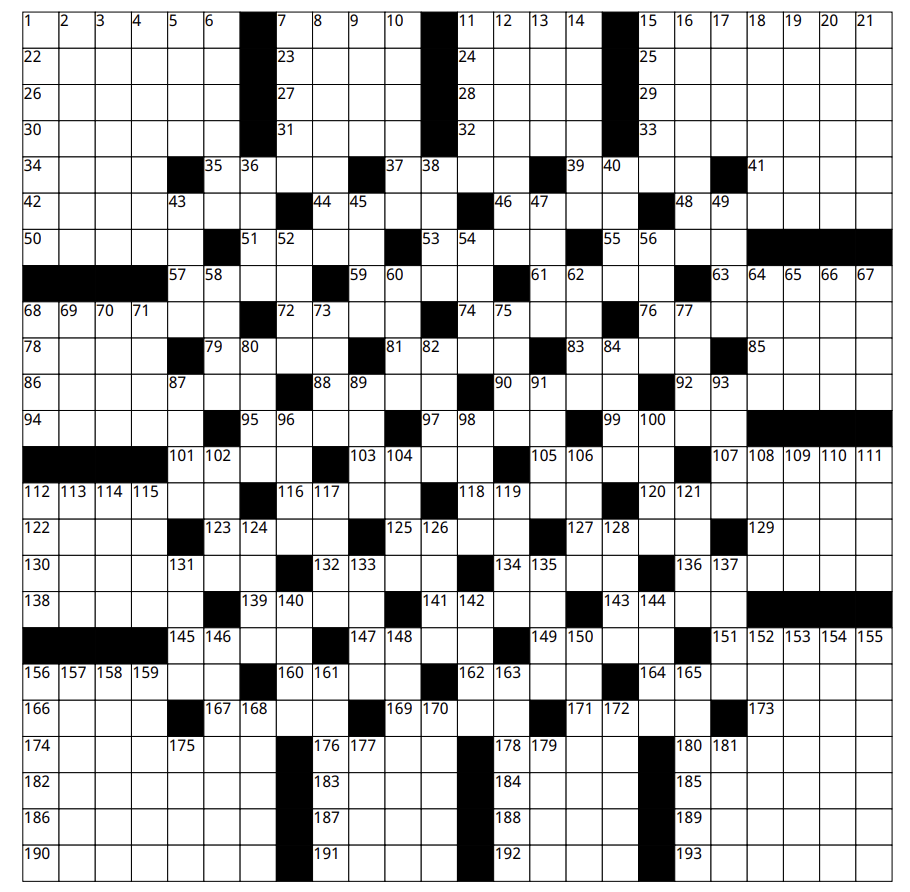}
  \caption{Example of $G(n)$ for $n\equiv4 \pmod{10}$}\label{24}
\endminipage\hfill
\minipage{0.48\textwidth}%
\centering
\caption*{$29\times29$}
  \includegraphics[width=.50\linewidth]{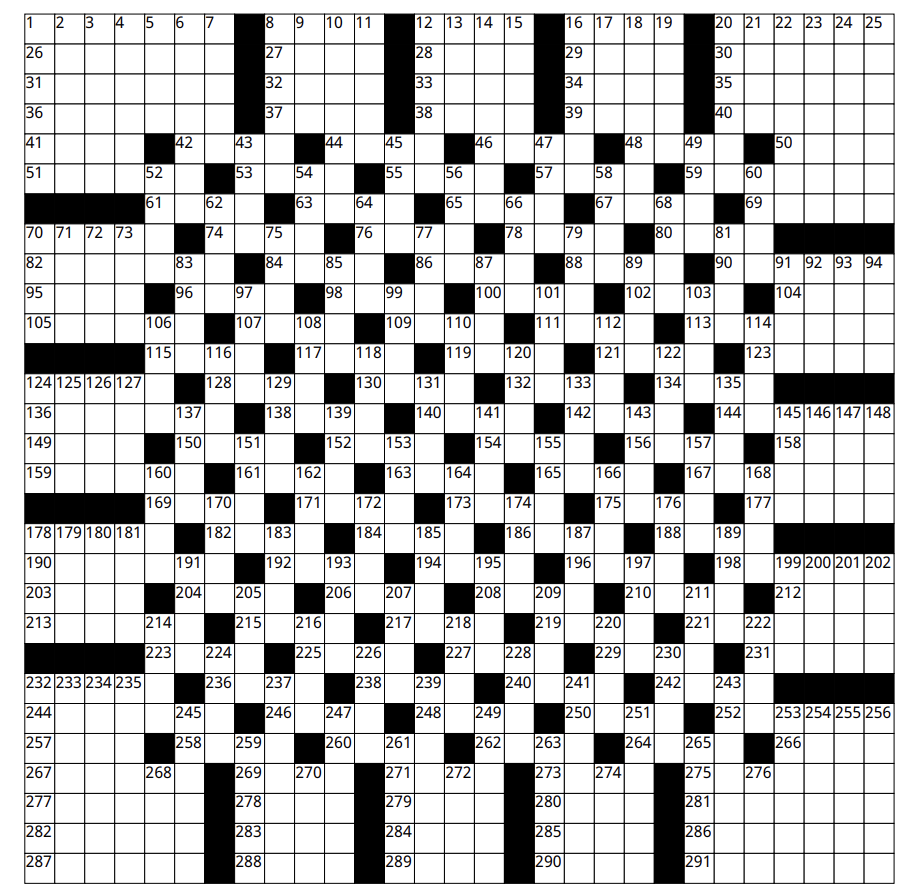}
  \caption{Example of $G(n)$ for $n\equiv9 \pmod{10}$}\label{29}
\endminipage
\end{figure}

We define $g(n)$ as the function to calculate the number of clues in $G(n).$

\section{Grid Calculations}
\label{g(n)}

For $n\geq 11$ we must break up the calculations of the grids into distinct cases depending on $n \pmod{5}$ to find the specific formulas for $G(n)$. 

A row or column is maximal if we have placed the maximum number of black squares in that row or column. Since we must satisfy Four$+$, the maximum number of clues in a row/column is $\lfloor \frac{n+1}{5}\rfloor $. We will show that in the construction of $G(n)$, every row/column is either maximal or one lower than maximal. The distinction in the cases $n \not \equiv 2 \pmod{5}$ is which rows/columns of the fundamental pane produce maximal rows/columns. Notice that the fundamental pane is 90-degree symmetric so the number of clues in each of the non-frame rows is the same as the number of clues in the corresponding non-frame column. For $n \equiv 2 \pmod{5}$, we also need to consider the middle strip.

First, to see why each row/column is either maximal or 1 less than maximal, note that each fundamental pane contributes one black square. We tile with $\lfloor \frac{n-8}{5}\rfloor$ fundamental panes in each row (and column by 90-degree symmetry), so the number of black squares $b_i$ in row (and column) $i$ is bounded by
$$
\lfloor\frac{n
+1}{5}\rfloor-2 \le \lfloor \frac{n-8}{5}\rfloor \le b_i \le \lfloor\frac{n
+1}{5}\rfloor-1
$$

\paragraph{$\mathbf{n\equiv 0 \pmod{5}}$} When $n \equiv 0 \pmod{5}$, we construct $\frac{n-10}{5}$ full fundamental panes. We also have one row tiled by the first row of the fundamental pane and one row tiled by the last row of the fundamental pane. If there is a black square in $(i,5)$ or $(i,6)$, then row $i$ will be maximal. Thus 2 of the 5 rows of the fundamental pane are maximal. The first and last row are submaximal. All 8 rows in the frame are maximal.

We then have $$g(n) =2\left( 8 \cdot \frac{n}{5} + \frac{n-10}{5} \cdot \left(2 \cdot\frac{n}{5} + 3\cdot \left(\frac{n}{5}-1\right)\right)+ 2 \cdot \left(\frac{n}{5}-1\right)\right)  = \frac{2}{5} \left(n^2-3n+20\right) \textrm{    if } n \equiv 0 \pmod{5}.$$

\paragraph{$\mathbf{n\equiv 1 \pmod{5}}$} Now we construct $\frac{n-11}{5}$ full fundamental panes plus an extra three rows: one maximal and two submaximal. There are 3 maximal rows in each fundamental pane. Each row of the frame has $\frac{n-1}{5}$ clues. Thus we see that
$$
g(n) = 2 \left(8 \cdot \frac{n-1}{5} + \frac{n-11}{5}\left(3 \cdot\frac{n-1}{5}  + 2 \cdot\left( \frac{n-1}{5}-1\right) \right) + 3 \cdot \frac{n-1}{5} -2 \right) = \frac{2}{5} \left(n^2-3n+12\right) \textrm{    if } n \equiv 1 \pmod{5}.
$$

\paragraph{$\mathbf{n\equiv 3 \pmod{5}}$} There are exactly $\frac{n-8}{5}$ full fundamental panes with no leftover rows. Each row in the pane is maximal as is each row in the frame. Thus 
$$
g(n) = 2 \cdot n \cdot \frac{n-3}{5} = \frac{2}{5}(n^2-3n) \textrm{    if } n \equiv 3 \pmod{5}.
$$

\paragraph{$\mathbf{n \equiv 4 \pmod{5}}$}
Each row in the frame contains $\frac{n-4}{5}$ clues. There are a total of $\frac{n-9}{5}$ fundamental panes, with 1 maximal row per frame. The one leftover row is maximal. Thus we have
$$
g(n) = 2 \left(8 \cdot \frac{n-4}{5} + \frac{n-9}{5} \left( \left(\frac{n-4}{5}+1\right) + 4\cdot \frac{n-4}{5}\right) + \frac{n-4}{5}+1 \right) = \frac{2}{5} (n^2-3n-4) \textrm{   if } n \equiv 4 \pmod{5}.
$$

\paragraph{$\mathbf{n \equiv 2 \pmod{5}}$} Finally, let's consider the last case, where we have enough additional rows after building fundamental panes to put a maximal strip in the middle. Note that this strip breaks 90-degree symmetry, so there will be a different number of row and column clues. Further, the construction of the strip depends on $n \pmod{10}$, so we get two different cases.

\paragraph{$\mathbf{n\equiv 2 \pmod{10}}$}
Each row in the frame is maximal. The rows in the $2 \times n$ strip are submaximal. There are a total of $\frac{n-12}{5}$ full fundamental panes containing 4 maximal rows, plus two leftover rows which are both maximal. Thus the number of clues in rows is
$$
8 \cdot \frac{n-2}{5} +\frac{n-12}{5} \cdot \left(4 \cdot \frac{n-2}{5} + \frac{n-2}{5}-1\right) + 2 \cdot \frac{n-2}{5} +2 \left(\frac{n-2}{5}-1\right)=\frac{n^2-3n+2}{5}.
$$

There are fewer clues in columns because of the middle strip. Now only 3 of the 5 columns in the fundamental pane are maximal. We get $\frac{n-7}{5}$ full fundamental panes across the rows plus 4 additional columns, 2 of which are maximal. Thus there are
$$
8\cdot \frac{n-2}{5} + \frac{n-12}{5} \left(3 \cdot \frac{n-2}{5} + 2\cdot \left(\frac{n-2}{5}-1\right)\right) + 4 \cdot \frac{n-2}{5} -2 =\frac{n^2-4n+14}{5}.
$$

Together we see that
$$
g(n) = \frac{1}{5} (2n^2-7n+16) \textrm{    if } n \equiv 2\pmod{10}.
$$

\paragraph{$\mathbf{n \equiv 7 \pmod{10}}$}
Notice that the column clue count for $n \equiv 2 \pmod{10}$ holds here as well. The difference in the row counts is that the middle $4 \times n$ strip contributes maximal rows. Since for $n \equiv 2\pmod{10}$, we had a $2\times n$ strip of submaximal rows and 2 additional maximal rows, our new construction has two additional maximal rows. Thus we can take the previous row count for $n \equiv 2 \pmod{10}$ and add 2. We get then that
$$
g(n) =\frac{1}{5} (2n^2-7n+16) + 2 = \frac{1}{5} (2n^2-7n + 21) \textrm{    if } n \equiv 7\pmod{10}.
$$

Combining all possible values of $n$ we obtain $g(n)$ to be 
\begin{equation*} g(n)= \begin{cases}
    \frac{2}{5}(n^2-3n+20), & n\equiv0 \pmod 5\\
     \frac{2}{5}(n^2-3n+12), & n\equiv1 \pmod 5\\
      \frac{1}{5}(2n^2-7n+14+(n\pmod{10}), & n\equiv2 \pmod 5\\
   \frac{2}{5}(n^2-3n), & n\equiv3 \pmod 5\\
   \frac{2}{5}(n^2-3n-4)  & n\equiv4 \pmod 5.\\
  \end{cases}
\end{equation*}

In each case we have verified Theorem \ref{THM:no14} and obtain the lower bound
\begin{equation*}
    g(n)\leq a(n).
\end{equation*}
 \begin{theorem}
  For each $n\geq11$, $G(n)$ is an $n\times n$ crossword grid satisfying Connectivity, Symmetry, and \fp { } that has $g(n)$ clues.\label{THM:no14}
    \end{theorem}

    We show all maximal grids for $n \le 30$ in Appendix \ref{APPE}. Some of the grids for $n < 11$ can be viewed as $G(n)$ constructions, though it is easier to count their clues and show maximality in other ways.
    
Note that for $n=14$ and $19$, our $G(n)$ construction is not maximal. Indeed, $g(14)=62<64=a(14)$ and $g(19) = 120 < 122 = a(19)$. For all $n \le 30$, we have verified $a(n)$ via linear programming using Claude Opus 5.5; see Isham's website for the code \cite{claudecode}. We have also included some combinatorial proofs for small $n$ that do not rely on computer proofs. See Appendices \ref{APPA} and \ref{APPB} for a proof of $a(10)$ and $a(14)$ respectively. We also are unable to prove $g(n) = a(n)$ when $n \equiv 0 \pmod{5}$. This is because our construction consists mostly of submaximal rows/columns and we do not have a method for proving why more maximal rows/columns are impossible.

\section{Optimizing $G(n)$}
There is a naive upper bound for $a(n)$ that satisfies \fp {} and does not depend on Connectivity or Symmetry. It follows from the fact that each clue takes at least 5 spaces -- four white spaces for the clue and one black space to separate it from the rest of the clues. This gives

\begin{equation} \label{eq:an2}
a(n)\leq 2n \left\lfloor\frac{n+1}{5} \right\rfloor= \begin{cases}
    \frac{2}{5}n^2, & n\equiv0 \pmod 5\\
   \frac{2}{5}(n^2-n), & n\equiv1 \pmod 5\\
   \frac{2}{5}(n^2-2n), & n\equiv2 \pmod 5\\
   \frac{2}{5}(n^2-3n), & n\equiv3 \pmod 5\\
    \frac{2}{5}(n^2+n), & n\equiv4 \pmod 5.\\
  \end{cases}
\end{equation}

We require the following lemma to show that in rows/columns 4 and $n-3$, there are either a maximal number of clues or one fewer than maximal. This helps us tighten the naive upper bound.

Observe that by \fp, if there is a black square in row 4, there must be a corresponding black square in rows 1, 2, and 3. Symmetrically, this is true for row $n-3$ and columns 4 and $n-3$. We will make use of this frequently in the following proofs.

\begin{lemma}\label{lem:submax}
   Let $n \ge 11$. Let $M_n$ be a $n\times n$ crossword with $a(n)$ clues. The number of clues in rows/columns 4 and $n-3$ of $M_n$ is either $\lfloor \frac{n+1}{5}\rfloor$ or $\lfloor \frac{n+1}{5}\rfloor-1$.
\end{lemma}

\begin{proof}
    Assume toward a contradiction that one of these rows/columns has $k< \lfloor \frac{n+1}{5}\rfloor-1$ clues. Suppose without loss of generality this occurs in row 4. Then by \fp, if there is a black square in $(i,4)$, there must also be one in $(i,3), (i,2),$ and $(i,1)$. Thus rows 1, 2, and 3 also have $k$ clues. In general, note that $r_1=r_2=r_3=r_4$ (as additional black squares, also known as ``cheater squares" do not contribute any additional clues). Note that adding one more black square to row 4 removes at most 1 clue from column 4, but adds one clue for rows 1, 2, 3, and 4. Thus we have increased the number of clues in $M_n$ by at least 3, contradicting that $M_n$ had the maximal number of clues.
\end{proof}

We can now prove Theorem \ref{THM:main}.
\begin{proof}[Proof of Theorem \ref{THM:main}]
 We start with the simplest case - when $g(n)$ is equal to the naive upper bound on $a(n)$.\\

\noindent\textbf{Case 1: $n \equiv3 \pmod{5}$}

   The naive upper bound in Equation \eqref{eq:an2} agrees with $g(n)$. Thus $a(n) = g(n)$.\\

\noindent \textbf{Case 2: $n \equiv1 \pmod{5}$}  For each black square in $(4,j)$ or $(n-3,j)$, observe that $c_j$ decreases by 1. Similarly, for each black square in $(i,4)$ or $(n-3,4)$, $r_i$ decreases by 1. We can assume that all rows and columns are maximal (so have $\frac{n-1}{5}$ clues) and subtract off the number of black squares in rows/columns 4 and $n-3$. 

By Lemma \ref{lem:submax}, these four columns are either maximal or 1 fewer than maximal. If row 4 is maximal, then it has $\frac{n-1}{5}-1$ black squares, so we lose $\frac{n-1}{5}-1$ maximal columns. If row 4 is submaximal, then it has $\frac{n-1}{5}-2$ black squares, so we lose $\frac{n-1}{5}-2$ maximal columns plus the one maximal row (row 4). Thus in either case we lose $\frac{n-1}{5}-1$ maximal rows/columns. The same logic holds for row $n-3$ and columns 4 and $n-3$. 

If we have a black square in both $(i,4)$ and $(i,n-3)$, then $r_i$ decreases by 2. So we can assume that we do not have both $(i,4)$ and $(i,n-3)$ black as the count is the same if they are both black. Symmetrically, we can assume the same for columns. Thus there are at most
$$
(2n) \cdot \frac{n-1}{5} - 4\left(\frac{n-1}{5}-1\right) = \frac{2n^2-2n-4n+24}{5} = \frac{2}{5}(n^2-3n+12)
$$
clues. \\

\noindent\textbf{Case 3: $n\equiv 4\pmod{5}$} Note that since $n\ne 14, 19$ and $n\ge 11$, then $n\ge 24$ in this case. Note that a row is maximal if and only if there are black squares in columns $c_5, c_{10}, \ldots, c_{n-4}$. There are no additional squares that can be black in this row. Symmetrically, a column is maximal if and only if there are black squares in rows $r_5, r_{10}, \ldots, r_{n-4}$.

Suppose we put black squares at \((4,5),(4,10),\ldots,(4,n-4)\). By Four+, we must also put black squares in the same positions in rows 1, 2, and 3. Thus we gain 4 clues in the rows, but force each of the \(\frac{n+1}{5}-1\) corresponding columns to lose at least one clue.

More generally, adding a maximal row in a position that was previously submaximal forces black squares into these same columns. Once a column has lost a clue, we may be able to add additional black squares without losing another clue. However, by Four+, we can have at most 4 consecutive black squares in this block before another black square forces the column to lose an additional clue.

Thus, a loss of one clue in each of the \(\frac{n+1}{5}-1\) columns can accommodate at most 4 additional maximal rows. These rows contribute at most 4 additional clues, while the corresponding columns lose at least \(\frac{n+1}{5}-1\) clues. Since \(n\ge24\),

\[
4-\left(\frac{n+1}{5}-1\right)\le0.
\]

Therefore, adding these maximal rows cannot increase the total number of clues. While the maximal rows do not need to be next to each other, this is clearly the optimal packing strategy - otherwise each black square added reduces the clues in the column by 1.

Thus we have at most $2 \cdot \left(\lfloor \frac{n+1}{5} \rfloor -1\right)$ maximal rows/columns. Therefore the number of clues in the crossword is at most
$$
2n \cdot \left(\lfloor \frac{n+1}{5}\rfloor-1\right) + 2 \cdot \left(\lfloor \frac{n+1}{5} \rfloor -1\right)
$$
assuming all $2n$ rows/columns are submaximal then adding the largest number of maximal rows/columns back. Simplifying, this becomes
$$
2n \frac{n-4}{5} + \frac{2(n-4)}{5} = \frac{2}{5} \left((n+1)(n-4)\right) = \frac{2}{5} (n^2-3n-4).
$$

\noindent\textbf{Case 4:} $n \equiv 2 \pmod{5}$.

Suppose that $n \equiv 2 \pmod{5}$. Observe that the maximum number of clues in a row or column is $\frac{n-2}{5}$. For each black square in row 4 (or $n-3$), observe that the corresponding column can only have at most $\frac{n-2}{5}-1$ clues. Symmetrically, for each black square in column 4 (or $n-3$), the corresponding row can only have at most $\frac{n-2}{5}-1$ clues. We can determine an upper bound on the number of clues by assuming that all $2n$ rows and columns are maximal, then subtract off the number of black squares in rows/columns 4 and $n-3$. Note that if these special rows and columns are maximal, then there are exactly $\frac{n-2}{5}-1$ black squares in each of these. However, note that by symmetry, if $n \equiv 7 \pmod{10}$, it is possible to line up the rows so that square $(4,j)$ and $(n-3,j)$ are both black, minimizing the number of submaximal rows. If $n \equiv 2 \pmod{10}$, we can do this for all but the middle two black squares. 

If any of these special four rows/columns is submaximal, then we can potentially add one more maximal column at the loss of a maximal row. Thus there is no gain to assuming these are submaximal.

Thus we have at most
$$
2n \cdot \frac{n-2}{5} - 3 \cdot \left(\frac{n-2}{5}-1\right) - f(n)
$$
clues, where $f(n) = 0$ if $n \equiv 7  \pmod{10}$ and $f(n) = 1$ if $n \equiv 2 \pmod{10}$. Simplifying, we have
$$
\frac{(2n^2-7n + 21) - 5f(n)}{5} = \frac{(2n^2-7n+14+ (n \hspace{-.05in}\pmod{10}))}{5}.
$$
Thus in this case $a(n) = g(n)$.

\end{proof}

\section{Grid Counting}

In \cite{ferland17}, Ferland describes ways to use symmetry to count the total number of Across and Down clues. We summarize his ideas here. Let $A_n$ be the total number of Across clues. Let $D_n$ be the total number of Down clues. Let $C_n$ be the total number of clues. 

We number the rows of the grid from top to bottom and the columns from left to right, beginning at 1. Note that $(i,j)$ would then refer to the square in row $i$ and column $j$. We denote the number of clues in row $i$ by $r_i$ and the number in the column $i$ by $c_i$. Furthermore, $n$ is the length/width of the grid.

 When $n$ is odd and our grid is a square, we denote the center row as $r_{\frac{n+1}{2}}$ and the center column as $c_{\frac{n+1}{2}}$. 
 When $n$ is even, we do not have a center row or column so instead we use the center line as a reference. From symmetry, we can count $A_n$ ($D_n$) using the first half of the rows (columns). The total number of clues is $C_n =A_n+D_n.$ 

 \begin{equation} 
    A_n=\begin{cases}
    2\displaystyle\sum\limits_{i=1}^{n/2} r_i, & \text{if \textit{n} is even}\\
    r_{(n+1)/2}+2\displaystyle\sum\limits_{i=1}^{(n-1)/2} r_i, & \text{if \textit{n} is odd}.
  \end{cases} 
  \label{2.2}
\end{equation}

\begin{equation*}
  D_n=\begin{cases}
    2\displaystyle\sum\limits_{i=1}^{n/2} c_i, & \text{if \textit{n} is even}\\
    c_{(n+1)/2}+2\displaystyle\sum\limits_{i=1}^{(n-1)/2} c_i, & \text{if \textit{n} is odd}.
  \end{cases}
\end{equation*}

Equation \eqref{2.2} leads directly to two results about the parity of clues. These results are originally found in a supplemental article \cite{supplemental} attached to \cite{ferland17} for $n=15$, but are presented more generally in \cite{article}.

\begin{lemma} \cite[Lemma 1]{article}
    Suppose n is odd and we have an n $\times$ n crossword grid satisfying the Symmetry rule. Then, $r_{(n+1)/2}$, the number of clues in the center row, and $c_{(n+1)/2}$, the number of clues in the center column, have the same parity .
    \label{L1}
\end{lemma}

\begin{corollary}\cite[Corollary 1]{article}
    The number of clues in any n $\times$ n crossword grid satisfying the Symmetry rule is always even.
    \label{C1}
\end{corollary}

\subsection{Counting the $9 \times 9$ Grid}

We see the utility of \eqref{2.2}, Lemma \ref{L1}, and Corollary \ref{C1} through counting the number of clues in a $9\times9$ crossword grid. This method of basic counting is inspired by Ferland's counting of the $15\times15$ grid \cite{supplemental}. \\

When we have a $9\times 9$ grid we utilize all the structure rules of American-Style crosswords detailed in Section \ref{Intro}. Note that as a consequence of \fp{} in a $9\times9$ grid, we say a ``TWO" is a row/column with exactly 2 clues. A ``ONE" is similarly defined.  

\begin{theorem}\label{thm:9}
    The maximum number of clues possible for a $9\times9$ crossword grid satisfying Connectivity, Symmetry, and \fp{ }is 26.
\end{theorem}

In the $9\times 9$ grid, $r_5$ is the center row and $c_5$ is the center column. By considering the cases in which the center square $(5, 5)$ is black or white, we get an immediate consequence of the Symmetry rule. Following Lemma \ref{L1}, we see that the number of clues in $r_5$ and $c_5$ have the same parity. Then by Corollary \ref{C1}, we know that since 9 is an odd number, a $9\times9$ grid will have an even number of clues. \\

We define the following as a result of \eqref{2.2}: 

\begin{align}\label{eq:3.1}
   & \text{Across}=r_{5}+2[r_1+r_2+r_3+r_4] \text{ and} \nonumber\\
    &\text{Down}=c_{5}+2[c_1+c_2+c_3+c_4.] 
\end{align}

\subsection{Proof of Theorem \ref{thm:9}}
\label{99proof section}
By Symmetry, we may assume that there are at least as many off-center TWO’s in the columns as the rows. We consider several cases, based on the placement of column TWO’s in the left-hand side of the grid. By \fp {} we cannot have 3 left TWO's.\\

\subsubsection{2 Left TWO's}
\begin{figure} [ht]
    \centering
    \includegraphics[width=0.25\linewidth]{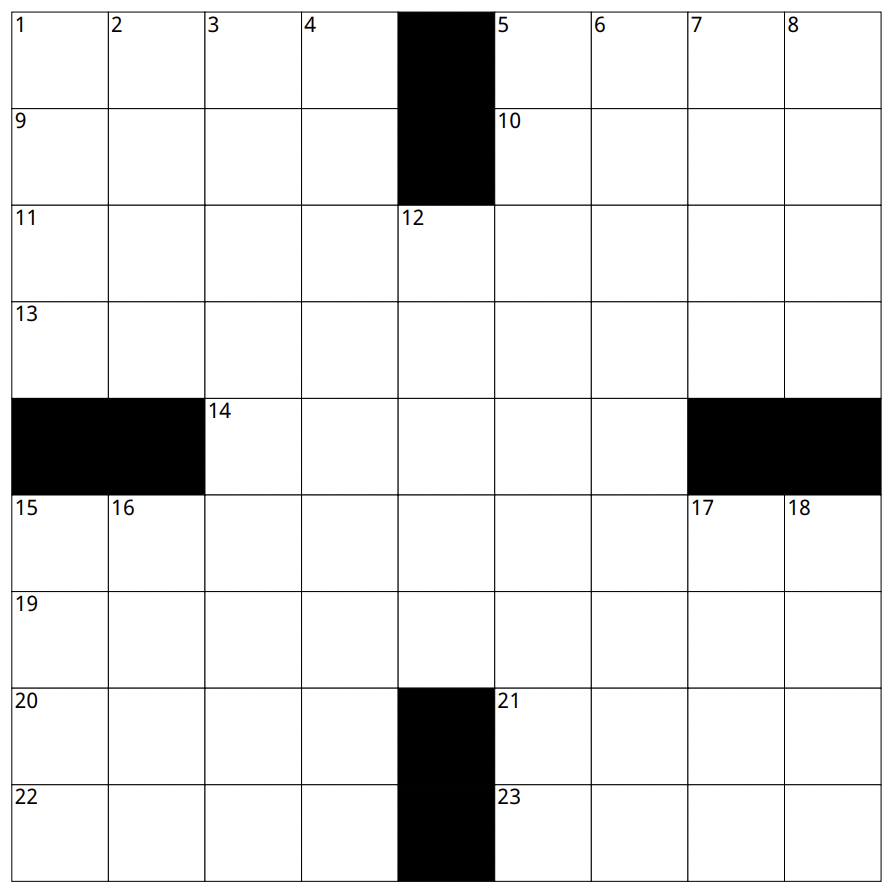}
    \caption{2 Left TWO's}
    \label{2Left}
\end{figure}

By Connectivity and \fp, the left TWO’s must be left-justified. We know that $r_5=1$ since the two left TWO's have to be in row 5 due to \fp. We also know that $c_5=1$ because the only way $c_5=2$ would be if $(5,5)$, the center square, was black which is impossible in the arrangement because it would violate \fp. Additionally, $r_1,r_2\leq2$ but all other rows and columns must equal 1 otherwise we would violate \fp{} and Connectivity. Figure \ref{2Left} shows such grid. 

Then by \eqref{eq:3.1} we have:

\begin{align*}
   & \text{Down}\leq 1+2[2(2)+2(1)]=13\\
    & \text{Across}\leq 1+2[2(2)+2(1)]=13\\
&      \Clues\leq26. \
\end{align*}

\subsubsection{1 Left TWO}
By Connectivity and \fp, the left TWO must be left-justified.

\begin{figure} [ht]
    \centering
    \includegraphics[width=0.25\linewidth]{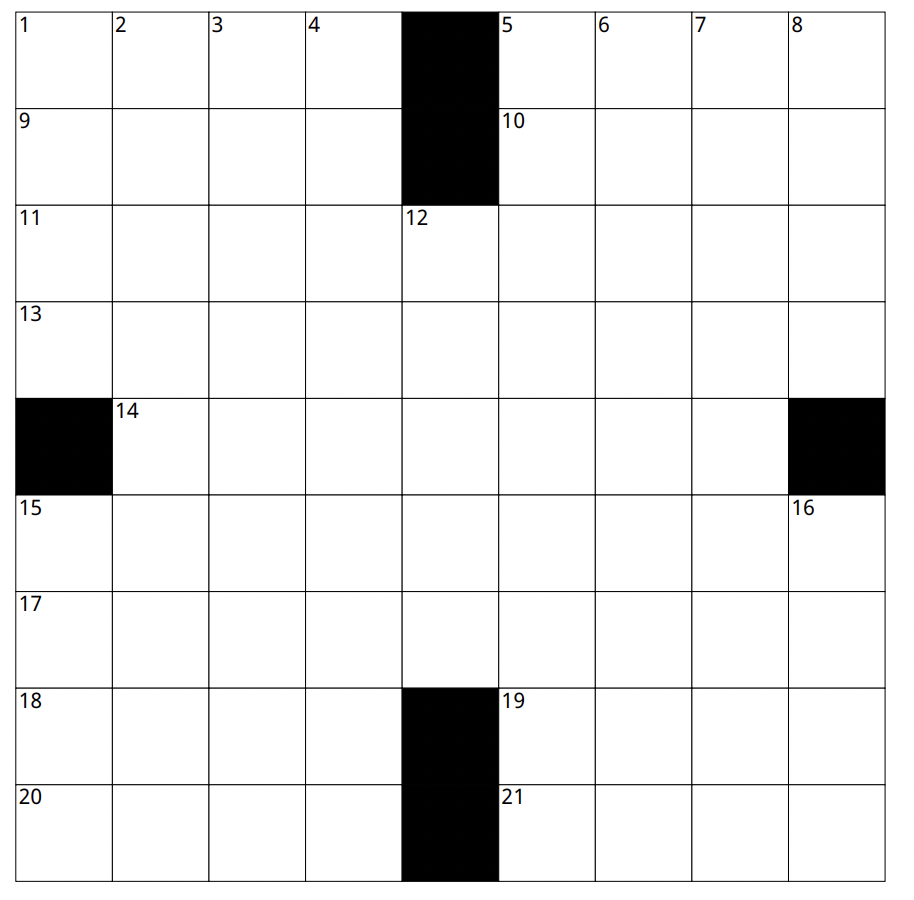}
    \caption{1 Left TWO}
    \label{1Left}
\end{figure}

By previous logic we know that $r_5,c_5=1$. Additionally, $r_1,r_2\leq2$ but all other rows and columns must equal 1 otherwise we would violate \fp{} and Connectivity. Figure \ref{1Left} shows such grid.  

Then by \eqref{eq:3.1} we have:

\begin{align*}
   & \text{Down}\leq 1+2[2+3(1)]=11\\
    & \text{Across}\leq 1+2[2(2)+2(1)]=13\\
&      \Clues\leq24. 
\end{align*}

\subsubsection{0 Left TWO's}
For 0 Left TWO's we refer to two cases: where $(5,5)$ is either black or white. \\

\begin{figure} [ht]
\centering
\subfloat[]{\includegraphics[width=0.21\linewidth]{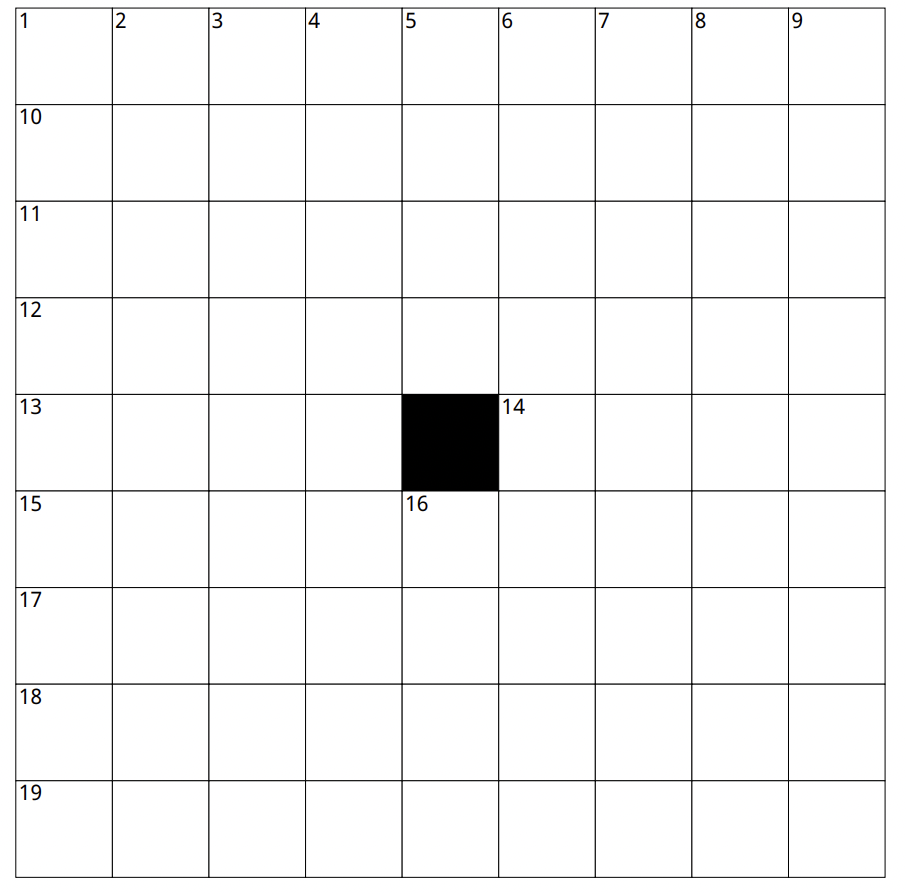}}\hspace{0.5cm}
\subfloat[]{  \includegraphics[width=0.21\linewidth]{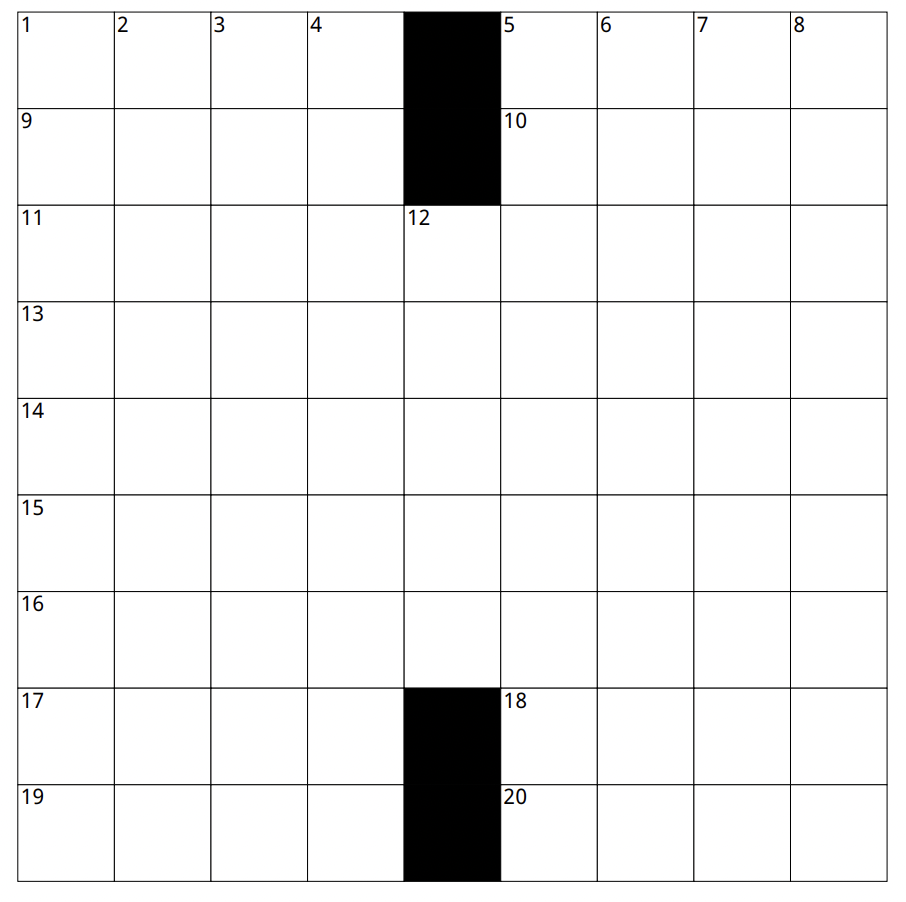}}\hspace{0.5cm}\\
\caption{0 Left TWO's} 
  \label{0l22}
\end{figure}

\underline{Case 1:} $(5,5)$ is black. Then by \fp, $c_5,r_5=2.$ We cannot place a single other black square on the grid without violating \fp. Figure \ref{0l22}a shows such grid.\\
Then by \eqref{eq:3.1} we have:

\begin{align*}
   & \text{Down}\leq 2+2[4(1)]=10\\
    & \text{Across}\leq 2+2[4(1)]=10\\
&      \Clues\leq20. \
\end{align*}

\underline{Case 2:} $(5,5)$ is white. Thus $r_5,c_5=1$. The only place to put a black square in the rows is in $r_1$ and $r_2$. So $r_1,r_2\leq2$. All other squares must be white. Figure \ref{0l22}b shows such grid.\\

Then by \eqref{eq:3.1} we have:

\begin{align*}
   & \text{Down}\leq 1+2[4(1)]=9\\
    & \text{Across}\leq 1+2[2(1)+2(2)]=13\\
&      \Clues\leq22. \
\end{align*}

Thus we have proven the maximum amount of clues in a $9\times9$ crossword puzzle satisfying \fp, Connectivity, and Symmetry is 26. The maximum amount of clues is achieved when we have 2 left TWO's because that is when there are the most amount of TWO's in the grid, the more TWO's there are means the higher the total clue count will be. Thus the result makes logistical sense.

As a result of this paper, the values of $a(n)$ up to $n=30$ have since been published by Abigail Schnitzer as sequence A378174 in Sloane’s On-Line Encyclopedia of Integer Sequences \cite{mysequence}. This is analogous to Pratt's sequence A243826 for American-style crosswords with Three+.

\section{Acknowledgements}
The authors thank Craig Larson and Bobby Jacobs for their encouragement and helpful conversations. This work is based on the Schnitzer's undergraduate thesis. The authors take responsibility for all results and writing. Claude Opus 5.5 was used to verify our results using ILP and Chat GPT 6 was used to check for typos and errors.

\begin{figure} 
    \centering
    \includegraphics[width=0.99\linewidth]{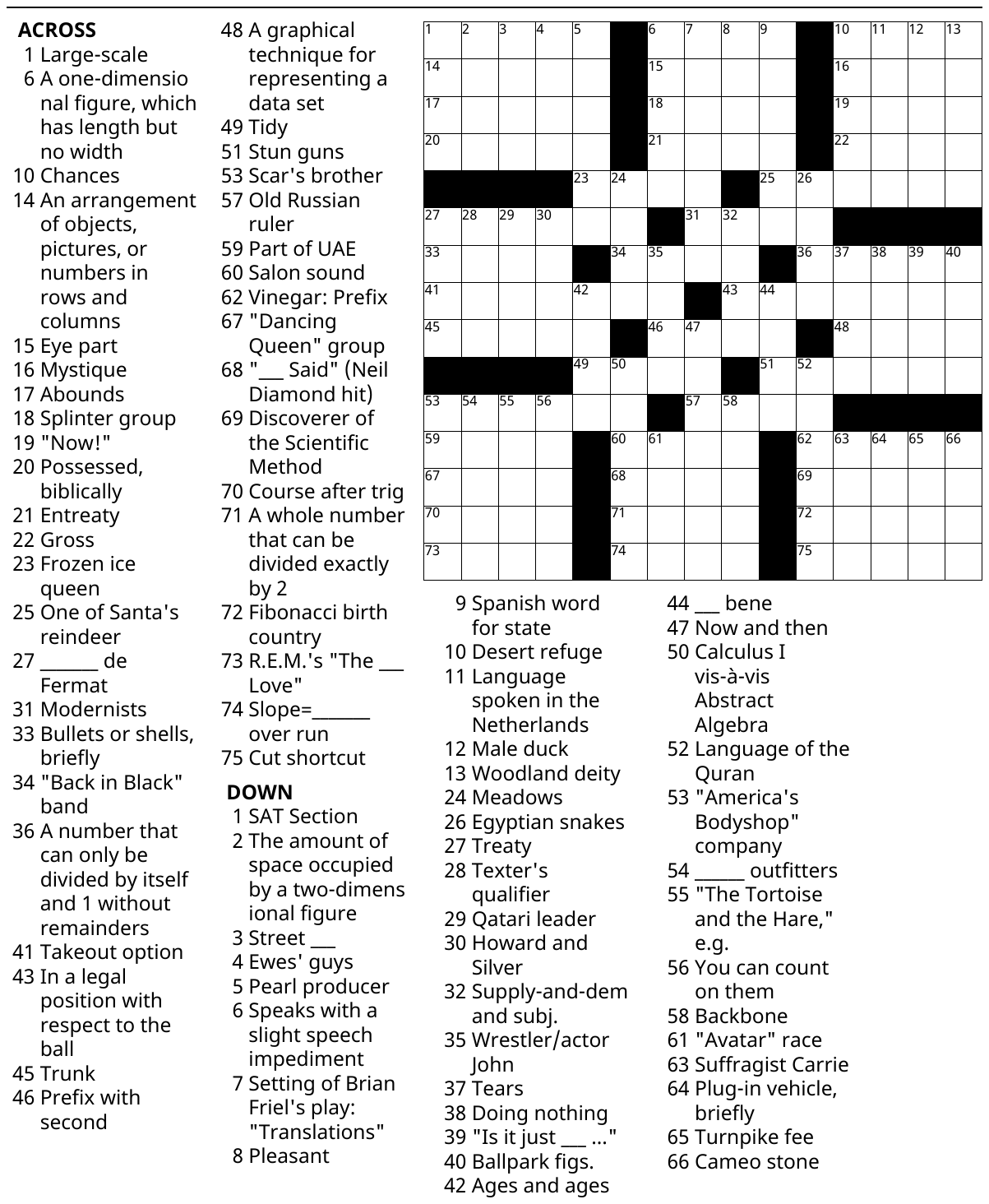}
    \caption{Maximal Crossword puzzle with 80 Clues. Solution appears in Appendix \ref{APPD}}
    \label{mathgrid}
\end{figure}

\newpage
\begin{center}
  \huge{\bf{Appendix}}
\end{center}
\appendix

 \section{$10\times10$ Maximal Grid}
\label{APPA}

\begin{figure}[ht]
\centering
 \includegraphics[width=0.25\linewidth]{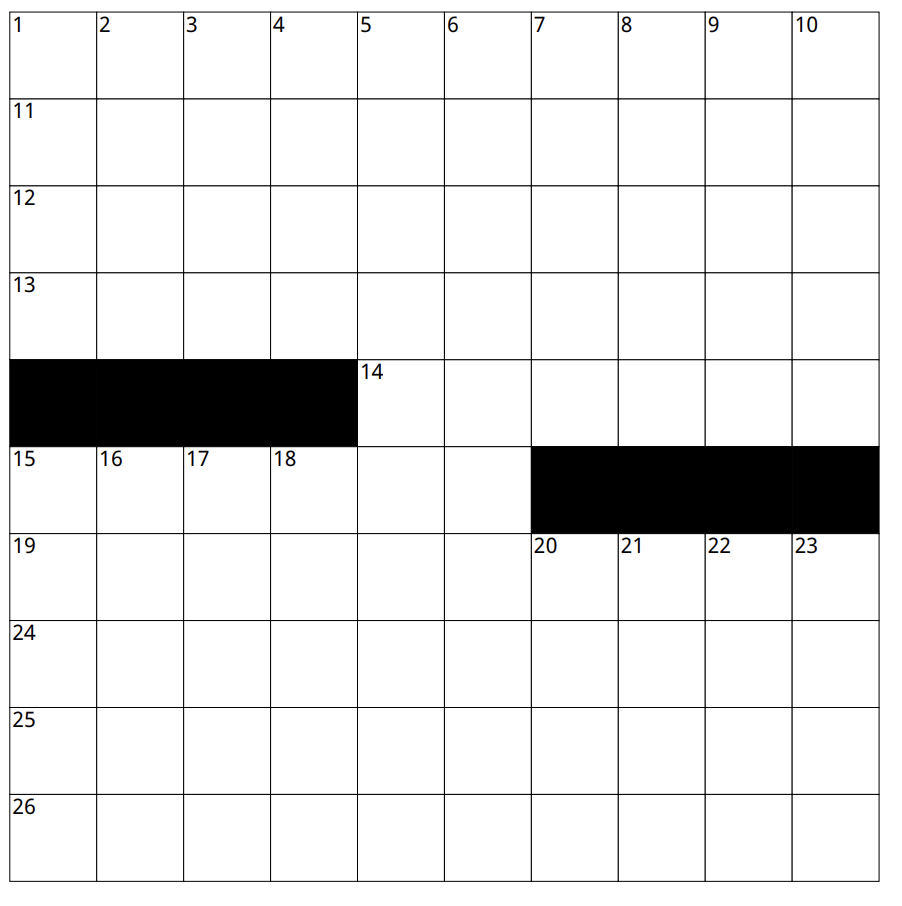}
 \caption{Assumed alignment}
    \label{2}
\end{figure}

For a $10\times10$ grid we employ the following equations as a result of the Symmetry rule: 

\begin{align} \label{eq:9}
   & \text{Across}=2[r_1+r_2+r_3+r_4+r_5] \text{ and} \nonumber\\
    &\text{Down}=2[c_1+c_2+c_3+c_4+c_5]\\
    &\Clues=\text{Across}+\text{Down}.\nonumber
\end{align}
By Symmetry, we may assume that there are at least as many off-center TWO’s in the columns as the rows. We consider several cases based on the placement of column TWO’s in the left-hand side of the grid. By Connectivity and \fp, we cannot have 5 left TWO's. It is important to note that we do not consider grids that are 180\textdegree { }rotations of each other as different grids. Thus we only calculate grids that are aligned like Figure \ref{2}. \\

\subsection{4 Left TWO'S}
By Connectivity and \fp, the 4 left TWO's must be left justified, as seen in Figure \ref{2}. We must have $r_5,c_5=1$ for the grid to satisfy Connectivity. We see that $r_1,r_2,r_3,r_4\leq2$.

Then by \eqref{eq:9} we have:

\begin{align*}
   & \text{Down}\leq 2[4(2)+1]=18\\
    & \text{Across}\leq 2[4(2)+1]=18\\
&      \Clues\leq36. \
\end{align*}
\\

\subsection{3 Left TWO'S}
\begin{figure} [ht]
    \centering
    \includegraphics[width=0.25\linewidth]{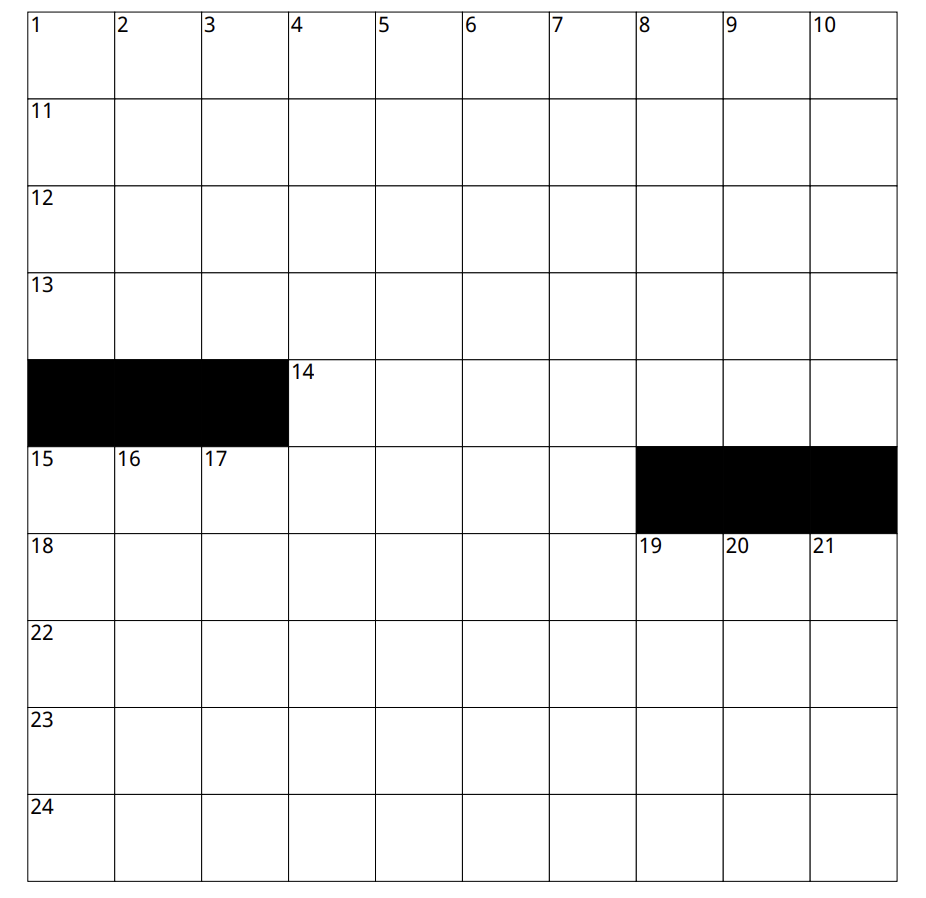}
    \caption{3 Left TWO's}
    \label{3l2}
\end{figure}

By Connectivity and \fp, the 3 left TWO's must be left justified, as seen in Figure \ref{3l2}. We must have $r_5,c_5=1$ for the grid to satisfy Connectivity. We see that $r_1,r_2,r_3,r_4\leq2$ and $c_4=1$. 

Then by \eqref{eq:9} we have:

\begin{align*}
   & \text{Down}\leq 2[3(2)+2]=16\\
    & \text{Across}\leq 2[4(2)+1]=18\\
&      \Clues\leq34. \
\end{align*}
\\

\subsection{2 Left TWO'S}
The 2 left TWO's need not be left justified and we refer to two possible cases. 

\begin{figure} [ht]
    \centering
    \includegraphics[width=0.25\linewidth]{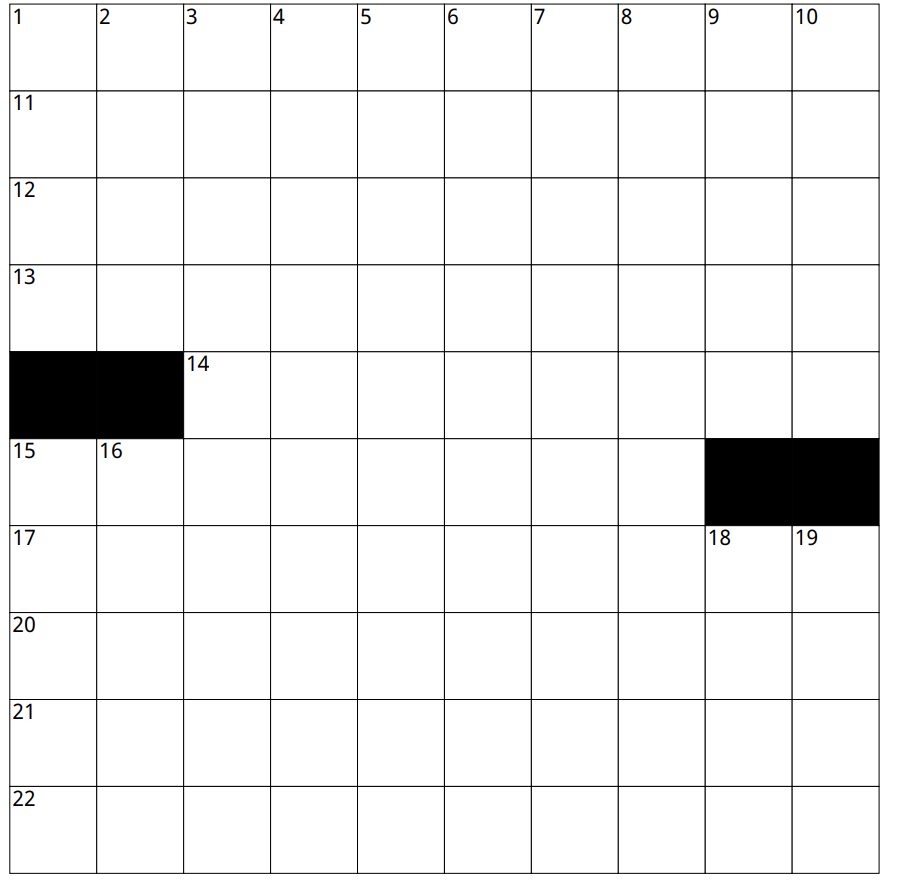}
    \caption{2 Left TWO's, Case 1}
    \label{2l2c1}
\end{figure}

\underline{Case 1:} The 2 Left TWO's are left justified, as seen in Figure \ref{2l2c1}. \\
 We must have $r_5,c_5=1$ for the grid to satisfy Connectivity. We see that $r_1,r_2,r_3,r_4\leq2$ and $c_3,c_4=1$. 

 Then by \eqref{eq:9} we have:

\begin{align*}
   & \text{Down}\leq 2[2(2)+3(1)]=14\\
    & \text{Across}\leq 2[4(2)+1]=18\\
&      \Clues\leq32. \
\end{align*}
\\

\begin{figure} [ht]
    \centering
    \includegraphics[width=0.25\linewidth]{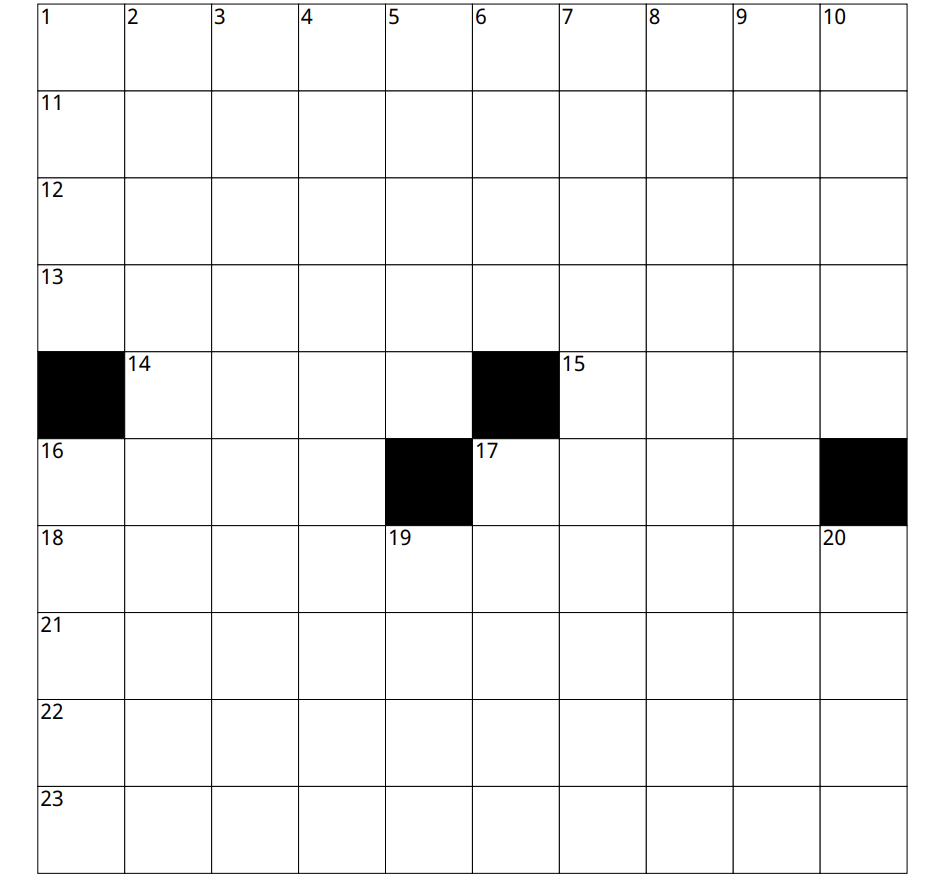}
    \caption{2 Left TWO's, Case 2}
    \label{2l2c2}
\end{figure}

\underline{Case 2:} The 2 left TWO's are not left justified, as seen in Figure \ref{2l2c2}. \\

The only space a black square can be without violating \fp {} is $(1,5)$. There can be no other squares that are black without violating Connectivity. 

Then by \eqref{eq:9} we have:

\begin{align*}
   & \text{Down}\leq 2[3(2)+2(1)]=16\\
    & \text{Across}\leq 2[2+4(1)]=12\\
&      \Clues\leq28. \
\end{align*}

\subsection{1 Left TWO'S}
The 1 left TWO need not be left justified and we refer to two possible cases.

\begin{figure} [ht]
    \centering
    \includegraphics[width=0.25\linewidth]{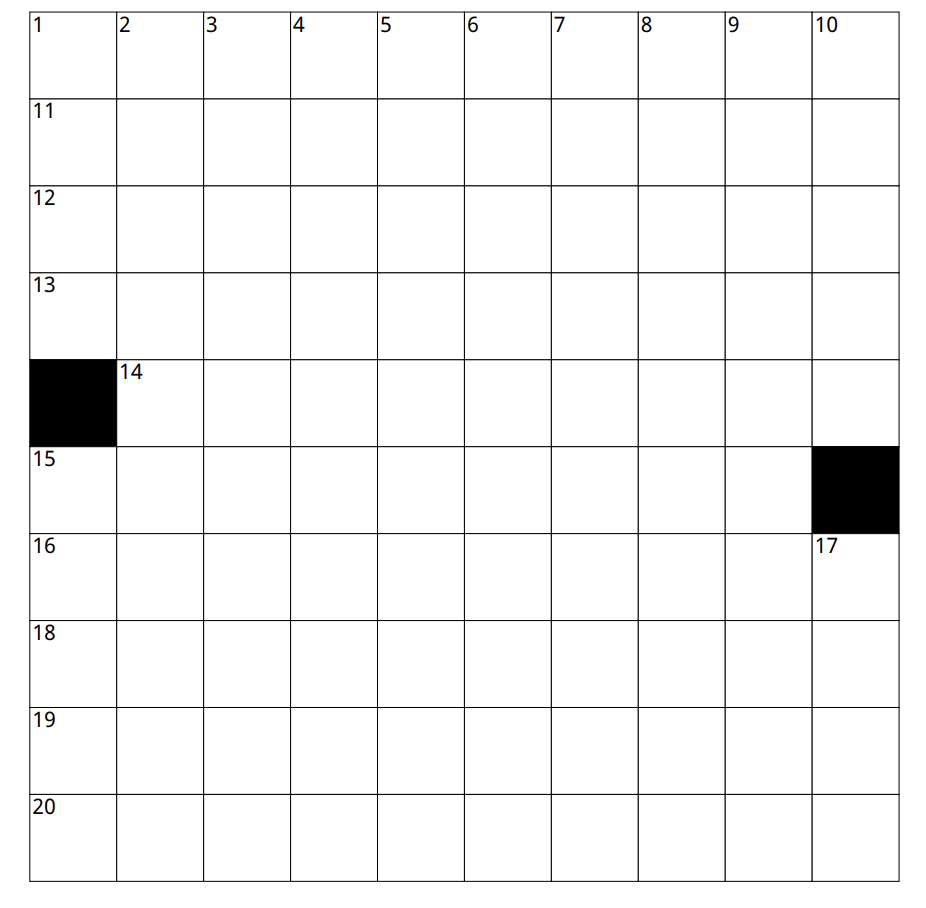}
    \caption{1 Left TWO, Case 1}
    \label{1l2c1}
\end{figure}

\underline{Case 1:} The left TWO is left justified, seen in Figure \ref{1l2c1}.

 We must have $r_5,c_5=1$ for the grid to have Connectivity. We see that $r_1,r_2,r_3,r_4\leq2$. All columns except for $c_1$ must $=1$ since we only have 1 left TWO. 

 Then by \eqref{eq:9} we have:

\begin{align*}
   & \text{Down}\leq 2[2+4(1)]=12\\
    & \text{Across}\leq 2[4(2)+1]=18\\
&      \Clues\leq30. \
\end{align*}
\\

\begin{figure}[ht]
    \centering
    \includegraphics[width=0.25\linewidth]{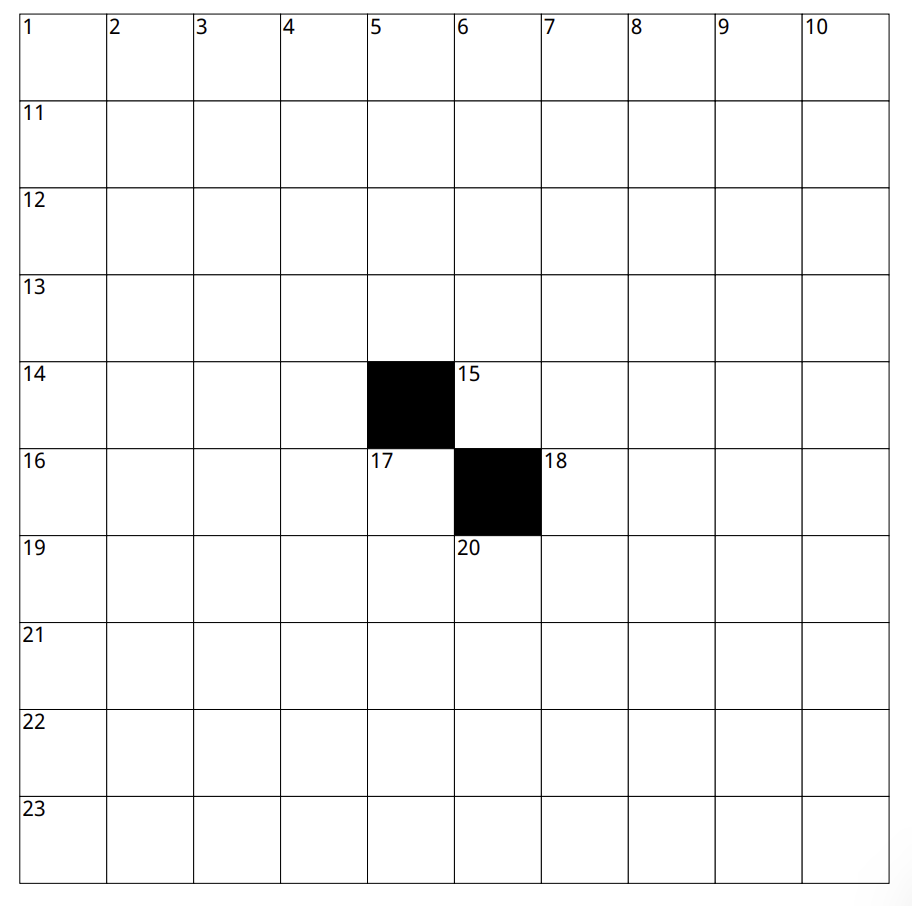}
    \caption{1 Left TWO, Case 2}
    \label{1l2c2}
\end{figure}

\underline{Case 2:} The 1 left TWO is not left justified, seen in Figure \ref{1l2c2}.

The only space a black square can be without violating \fp {} is $(1,6)$. There can be no other squares that are black without violating Connectivity. 

 Then by \eqref{eq:9} we have:

\begin{align*}
   & \text{Down}\leq 2[2(2)+3(1)]=14\\
    & \text{Across}\leq 2[2+4(1)]=12\\
&      \Clues\leq26. \
\end{align*}

\subsection{0 Left TWO'S}
\begin{figure} [ht]
    \centering
    \includegraphics[width=0.25\linewidth]{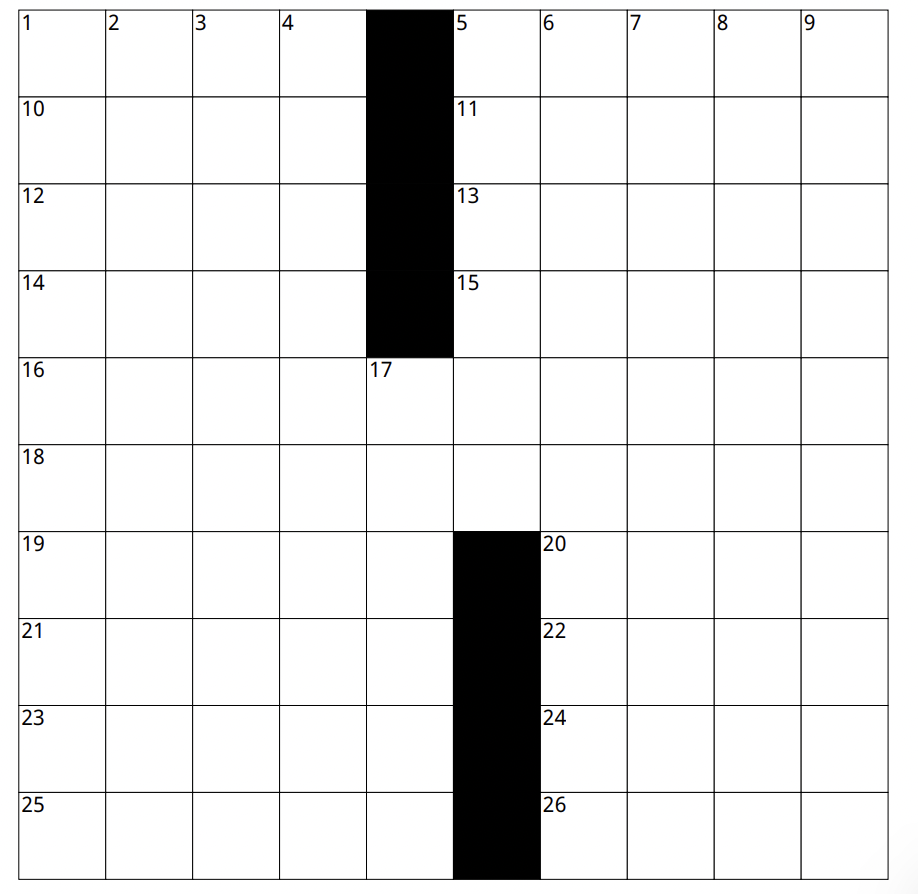}
    \caption{0 Left TWO's}
    \label{0l2}
\end{figure}

The only place back squares can be without violating the 0 left TWO's arrangement, Connectivity or \fp{} is $(1,5), (2,5), (3,5), (4,5)$, as seen in Figure \ref{0l2}. This yields that $r_1,r_2,r_3,r_4\leq 2$ and $r_5=1$. 

 Then by \eqref{eq:9} we have:

\begin{align*}
   & \text{Down}\leq 2[5(1)]=10\\
    & \text{Across}\leq 2[4(2)+1]=18\\
&      \Clues\leq28. \
\end{align*}
\\

Thus  we have proven the maximum amount of clues in a $10\times10$ crossword puzzle satisfying \fp, Connectivity, and Symmetry is 36.

\section{$14\times14$ Maximal Grid}
\label{APPB}

For a $14\times14$ grid we employ the following equations as a result of the Symmetry rule: 

\begin{align} \label{eq:14}
   & \text{Across}=2[r_1+r_2+r_3+r_4+r_5+r_6+r_7] \text{ and} \nonumber\\
    &\text{Down}=2[c_1+c_2+c_3+c_4+c_5+c_6+c_7]\\
    &\Clues=\text{Across}+\text{Down}.
    \nonumber    
\end{align}

By Symmetry, we may assume that there are at least as many off-center THREE’s in the columns as the rows. We consider several cases, based on the placement of column THREE’s in the left-hand side of the grid. By Connectivity and \fp, we cannot have 6 left THREE's.

\begin{figure} [ht]
    \centering
    \includegraphics[width=0.25\linewidth]{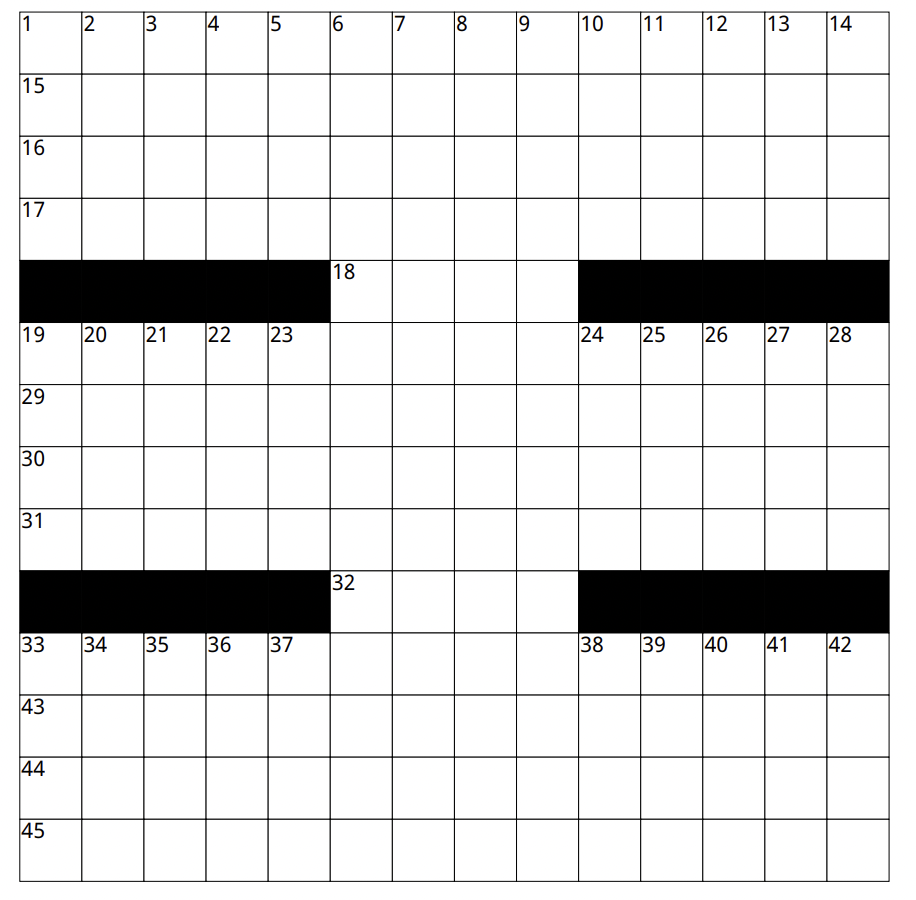}
    \caption{5 Left THREE's}
    \label{5l3}
\end{figure}
\subsection{5 Left THREE'S}
By \fp, the 5 left THREE's must be left justified, as seen in Figure \ref{5l3}. Thus, $r_5=1$. We see that $r_1,r_2,r_3,r_4\leq2$. In order to preserve \fp and connectivity, both $r_6,r_7 \le 2$. Notice that $c_i \le 3$ for $i=1, \ldots, 5$ but $c_6, c_7 \le 2$.

Then by \eqref{eq:14} we have:

\begin{align*}
   & \text{Down}\leq 2[2(2)+5(3)]=38\\
    & \text{Across}\leq 2[6(2)+1]=26\\
&      \Clues\leq64. \
\end{align*}

\begin{figure} [ht]
    \centering
    \includegraphics[width=0.25\linewidth]{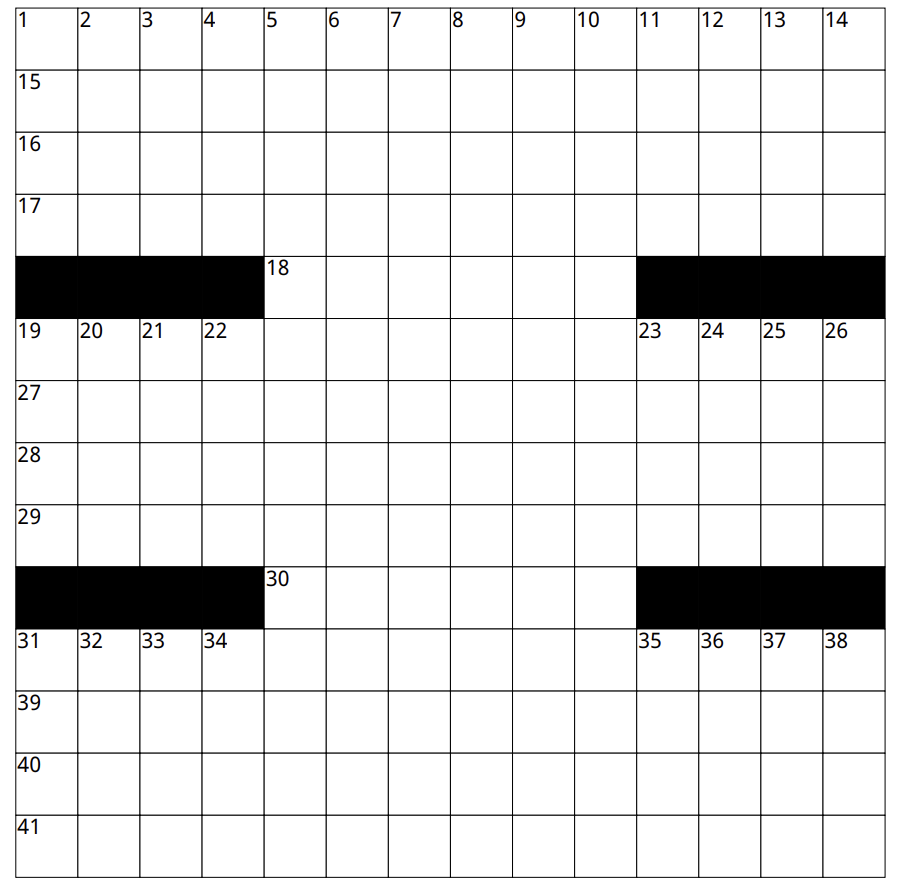}
    \caption{4 Left THREE's}
    \label{4l3}
\end{figure}
\subsection{4 Left THREE'S}
By \fp, the 4 left THREE's must be left justified, as seen in Figure \ref{4l3}. Thus, $r_5=1$. If $r_4=2$, then $r_1,r_2,r_3=2$ as well. That would yield $r_1+r_2+r_3+r_4=8$. However, if $r_4=1,$ then $r_1,r_2,r_3\leq3$. This would yield of maximum of $r_1+r_2+r_3+r_4\leq10$. Since $10>8$, we will assume, without loss of generality, that $r_4=1$. There is no way to make $r_6,r_7=3$ without violating \fp{} so $r_6,r_7\leq2$. With $r_1,r_2,r_3\leq3$, $c_5=1$.

Then by \eqref{eq:14} we have:

\begin{align*}
   & \text{Down}\leq 2[2(2)+4(3)+1]=34\\
    & \text{Across}\leq 2[3(3)+2(1)+2(2)]=30\\
&      \Clues\leq64. \
\end{align*}

\subsection{3 Left THREE'S}
\begin{figure} [ht]
\centering
\subfloat[]{\includegraphics[width=0.23\linewidth]{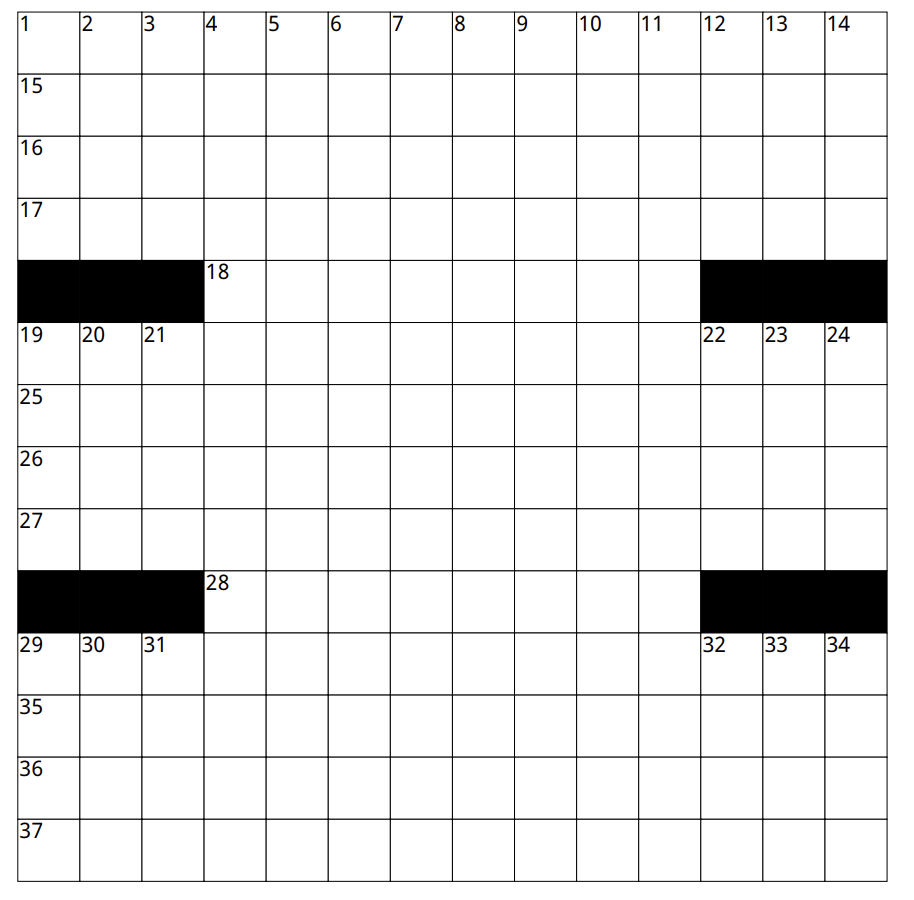}}\hspace{0.5cm}
\subfloat[]{ \includegraphics[width=0.23\linewidth]{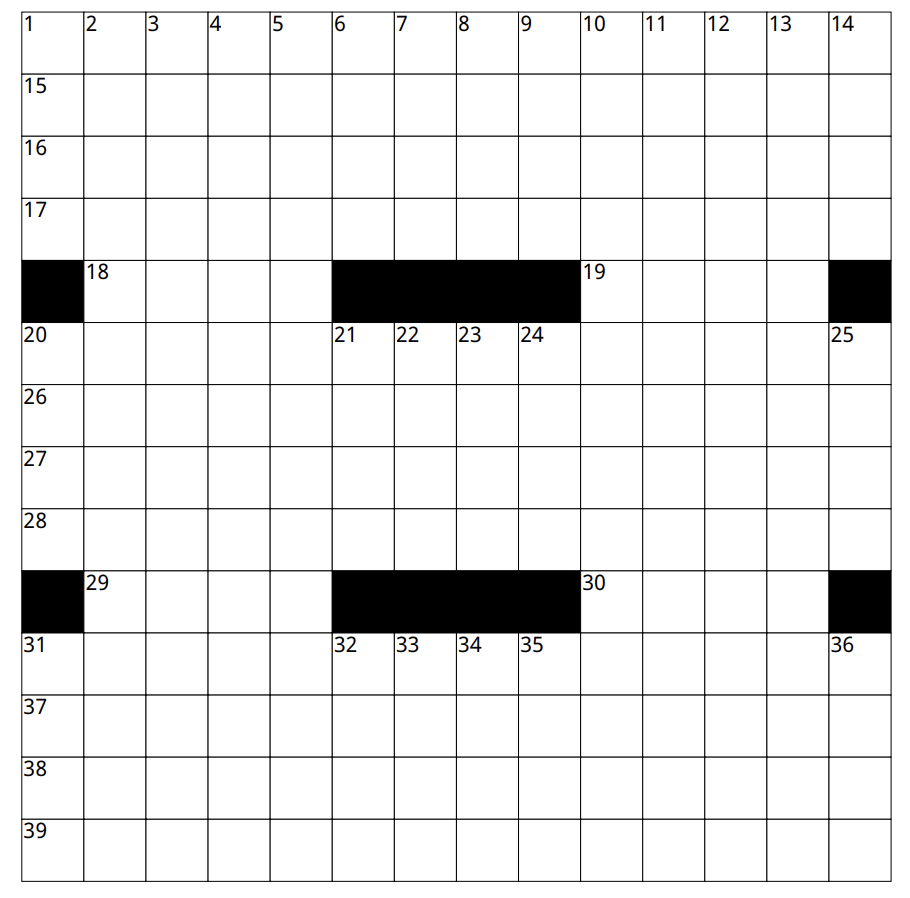}}\\
\subfloat[]{  \includegraphics[width=0.23\linewidth]{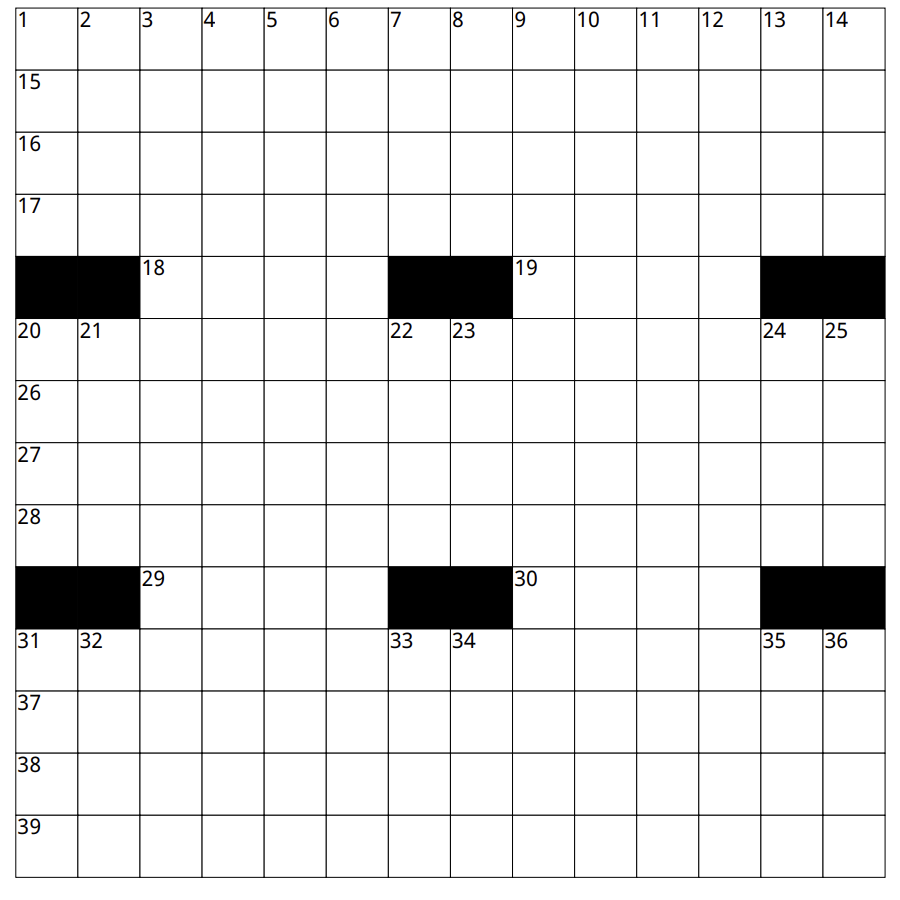}}\hspace{0.5cm}\\
\caption{3 Left THREE's} 
  \label{3l3}
\end{figure}

There are three possible cases of such arrangement, as seen in Figure \ref{3l3}. 

\underline{Case 1:} The 3 left THREE's are left justified, as seen in Figure \ref{3l3}a. By \fp, $r_5,c_4=1$. We see that $r_1,r_2,r_3,r_4\leq3$. However when $r_1,r_2,r_3,r_4=3$, then $c_5=1$.

Then by \eqref{eq:14} we have:

\begin{align*}
   & \text{Down}\leq 2[2(2)+3(3)+2(1)]=30\\
    & \text{Across}\leq 2[4(3)+1+2(2)]=34\\
&      \Clues\leq64. \
\end{align*}
\\

\underline{Case 2:} One left THREE is left justified, as seen in Figure \ref{3l3}b. Clearly, $r_5=2$. If $r_6\leq3$, then $r_7\leq2$ and vice versa. In order to maximize, we assume that $r_1\leq3$. No other squares can be black without violating Connectivity or \fp. 

Then by \eqref{eq:14} we have:

\begin{align*}
   & \text{Down}\leq 2[3(3)+3(1)+2]=28\\
    & \text{Across}\leq 2[2(2)+2(3)+3(1)]=26\\
&      \Clues\leq54. \
\end{align*}
\\

\underline{Case 3:} Two left THREE's are left justified, as seen in Figure \ref{3l3}c. Clearly $r_5=1$. Then $r_1,r_2,r_3,r_4\leq3$. If $r_6\leq3$ then $r_7\leq2$ and vice versa.  When $r_1,r_2,r_3,r_4=3$, $c_5=1$. 

Then by \eqref{eq:14} we have:

\begin{align*}
   & \text{Down}\leq 2[3(3)+3(1)+2]=28\\
    & \text{Across}\leq 2[4(3)+3(2)]=36\\
&      \Clues\leq64. \
\end{align*}
\\

\subsection{2 Left THREE'S}
\begin{figure} [ht]
\centering
\subfloat[]{ \includegraphics[width=0.25\linewidth]{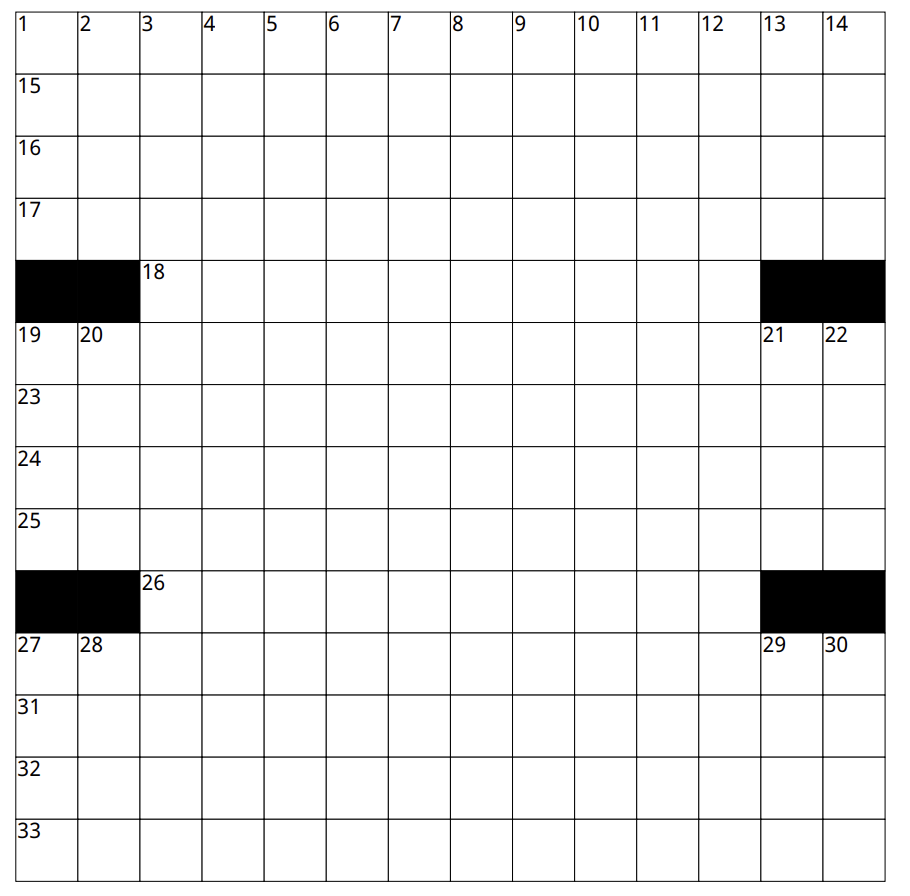}}\hspace{0.5cm}
\subfloat[]{  \includegraphics[width=0.25\linewidth]{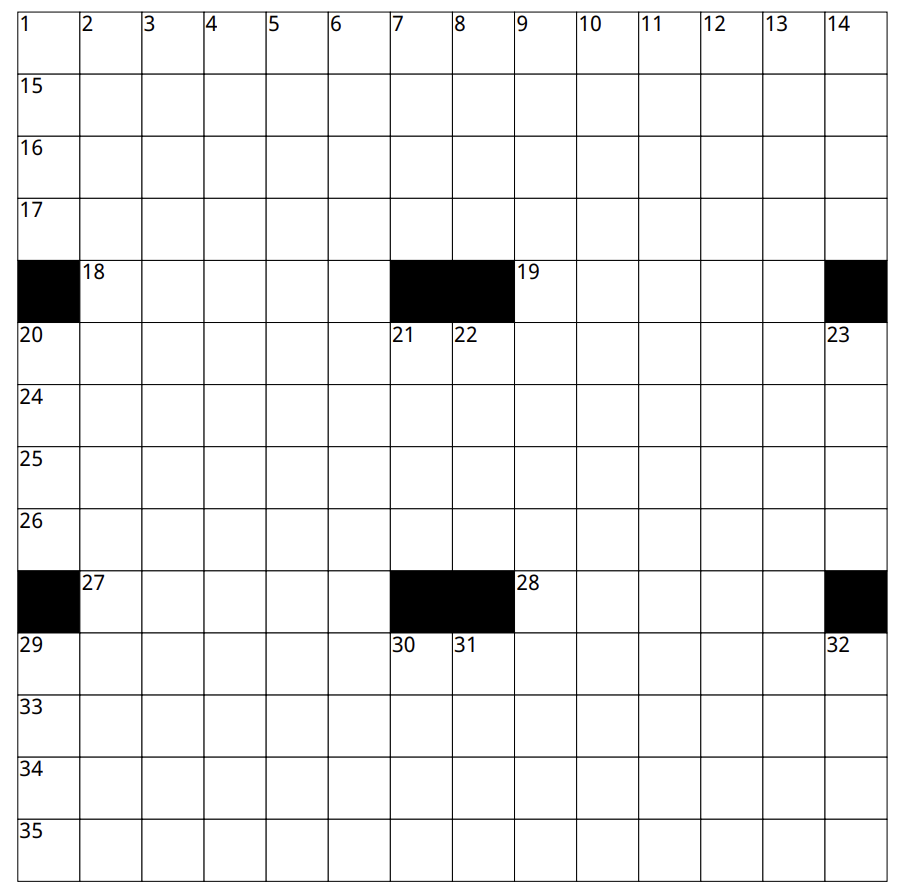}}\\
\subfloat[]{ \includegraphics[width=0.25\linewidth]{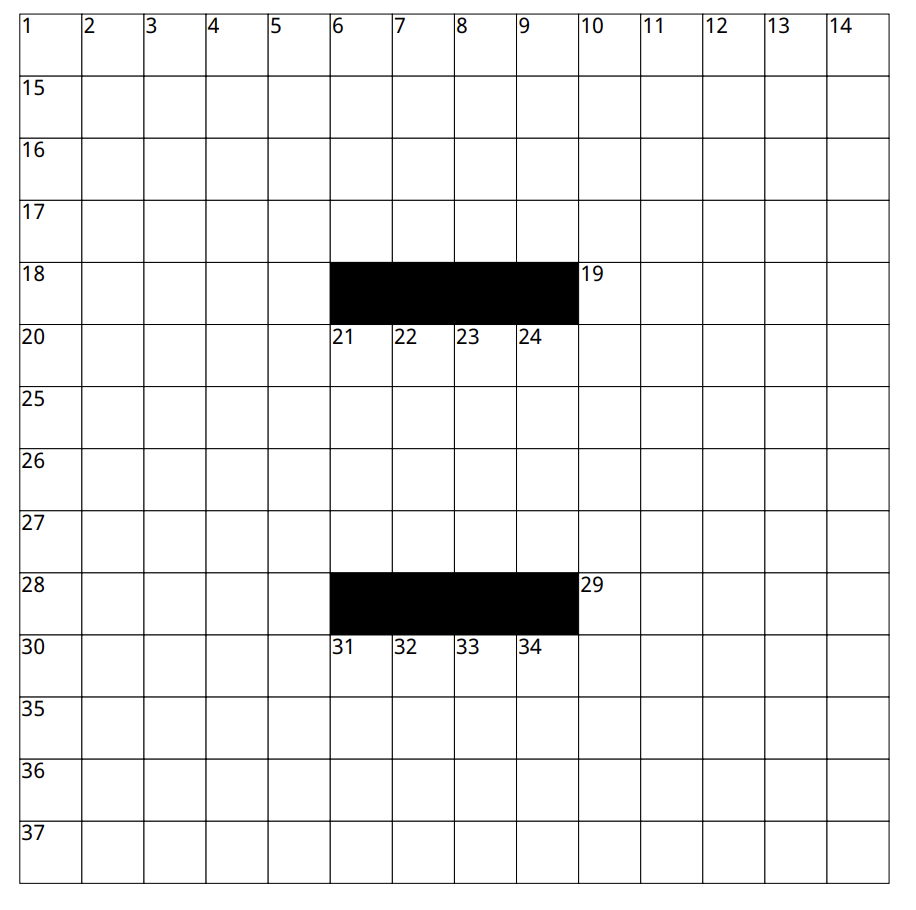}}\hspace{0.5cm}\\
\caption{2 Left THREE's} 
  \label{2l3}
\end{figure}

We again refer to 3 possible cases, as seen in Figure \ref{2l3}.

\underline{Case 1:} The 2 left THREE's are left justified, as seen in Figure \ref{2l3}a.  Then $r_1,r_2,r_3,r_4\leq3$. If $r_1,r_2,r_3,r_4=3$, then $c_5=1$, and this is clearly maximal since we want as many rows/columns to be THREE's. 

Then by \eqref{eq:14} we have:

\begin{align*}
   & \text{Down}\leq 2[2(3)+3(1)+2(2)]=26\\
    & \text{Across}\leq 2[4(3)+3(2)]=36\\
&      \Clues\leq62. \
\end{align*}
\\

\underline{Case 2:} One left THREE is left justified, as seen in Figure \ref{2l3}b.  Then $r_1,r_2,r_3,r_4\leq3$. If $r_1,r_2,r_3,r_4=3$, then $c_5=1$, and this is clearly maximal since we want as many rows/columns to be THREE's. 

Then by \eqref{eq:14} we have:

\begin{align*}
   & \text{Down}\leq 2[2(3)+4(1)+2]=24\\
    & \text{Across}\leq 2[4(3)+3(2)]=36\\
&      \Clues\leq60. \
\end{align*}
\\

\underline{Case 3:} No left THREE's are left justified, as seen in Figure \ref{2l3}c. Now only $r_1,r_2,r_3\leq3$ because $r_4\leq2$ in order to preserve Connectivity. Still $c_5=1$. In order to preserve Connectivity and \fp, $c_1,c_2,c_3,c_4\leq2$. 

Then by \eqref{eq:14} we have:

\begin{align*}
   & \text{Down}\leq 2[4(2)+2(3)+1]=26\\
    & \text{Across}\leq 2[3(3)+3(2)+1]=32\\
&      \Clues\leq58. \
\end{align*}

\subsection{1 Left THREE}
\begin{figure} [ht]
\centering
\subfloat[]{  \includegraphics[width=0.25\linewidth]{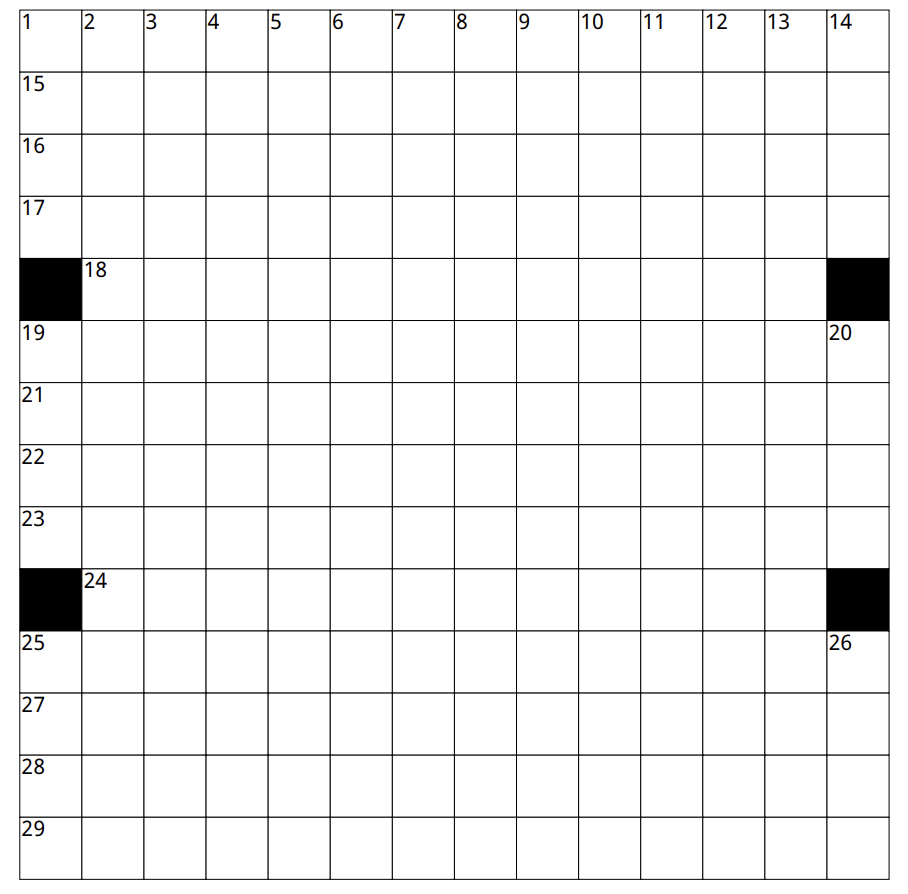}}\hspace{0.5cm}
\subfloat[]{  \includegraphics[width=0.25\linewidth]{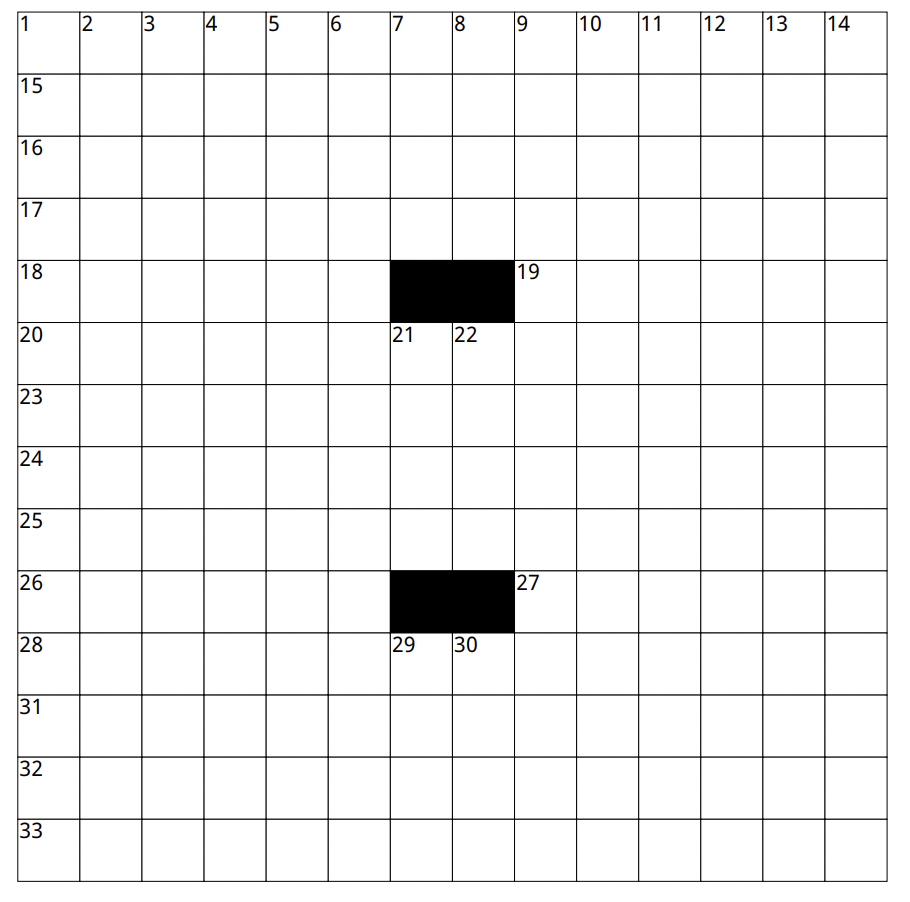}}\hspace{0.5cm}\\
\caption{1 Left THREE} 
  \label{1l3}
\end{figure}

We refer to two possible cases, as seen in Figure \ref{1l3}.

\underline{Case 1:} The 1 left THREE is left justified, as seen in Figure \ref{1l3}a. We see that $r_1,r_2,r_3,r_4\leq3$. If $r_1,r_2,r_3,r_4=3$, then $c_5=1$. There is no place to put a black square to split up $c_2,c_3,c_4$ without violating \fp, so $c_2,c_3,c_4=1$. Also, $c_6,c_7\leq2.$

Then by \eqref{eq:14} we have:

\begin{align*}
   & \text{Down}\leq 2[3+4(1)+2(2)]=22\\
    & \text{Across}\leq 2[4(3)+3(2)]=36\\
&      \Clues\leq58. \
\end{align*}
\\

\underline{Case 2:} The 1 left THREE is not left justified, as seen in Figure \ref{1l3}b. In order to maximize the amount of clues, $r_1,r_2,r_3,r_4\leq3$ and $c_1,c_2,c_3,c_4\leq2$. This means that $c_5=1$. 

Then by \eqref{eq:14} we have:

\begin{align*}
   & \text{Down}\leq 2[3+5(2)+1]=28\\
    & \text{Across}\leq 2[4(3)+2(2)+1]=34\\
&      \Clues\leq62. \
\end{align*}

\subsection{0 Left THREE'S}
We see that $r_1,r_2,r_3,r_4,r_5\leq3$. Still, $c_1,c_2,c_3,c_4\leq2$. This means that $c_5=1$ and $c_6,c_7\leq2$ as well. 

Then by \eqref{eq:14} we have:

\begin{align*}
   & \text{Down}\leq 2[6(2)+1]=26\\
    & \text{Across}\leq 2[5(3)+2+1]=36\\
&      \Clues\leq62. \
\end{align*}
\\

Thus the maximum number of clues in a $14\times14$ crossword grid satisfying Connectivity, Symmetry, and \fp, is $64$.

\clearpage
\section{Solution to Figure \ref{mathgrid}}
\label{APPD}

\begin{figure} [ht]
    \centering
    \includegraphics[width=0.7\linewidth]{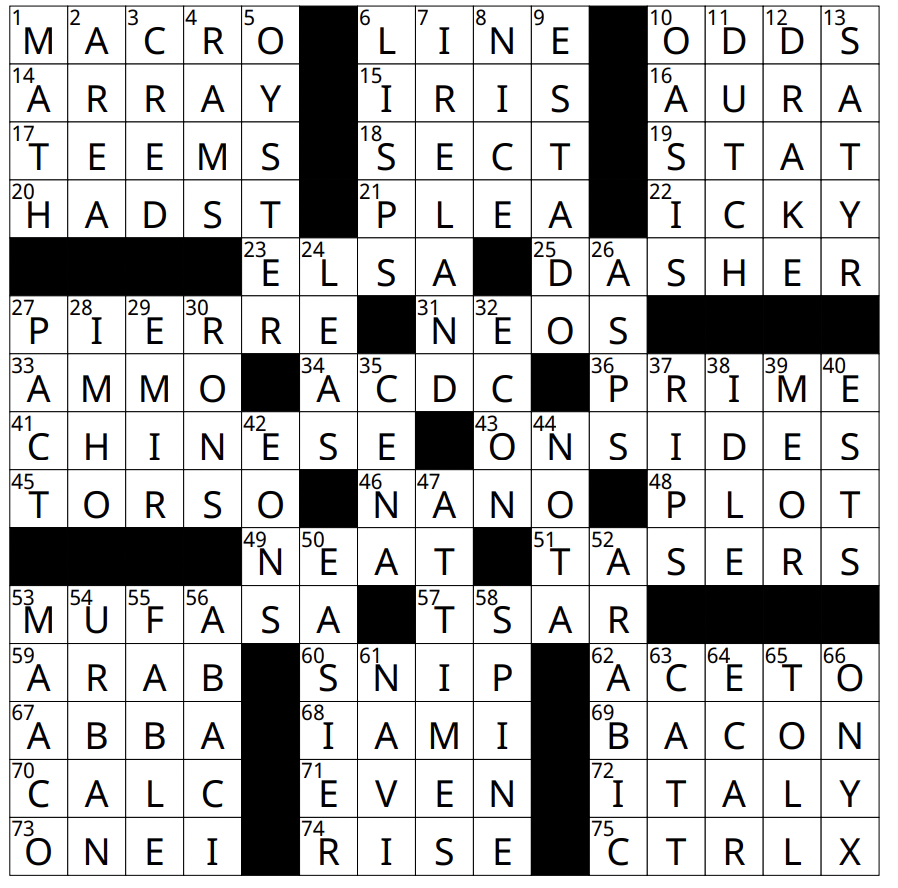}
    \caption{Solution to Figure \ref{mathgrid}}
    
\end{figure}
\clearpage
\section{Maximal Grids $[4,30]$}
\label{APPE}

\begin{figure}[!htb]
   \begin{minipage}{0.48\textwidth}
     \centering
      \includegraphics[width=0.35\linewidth]{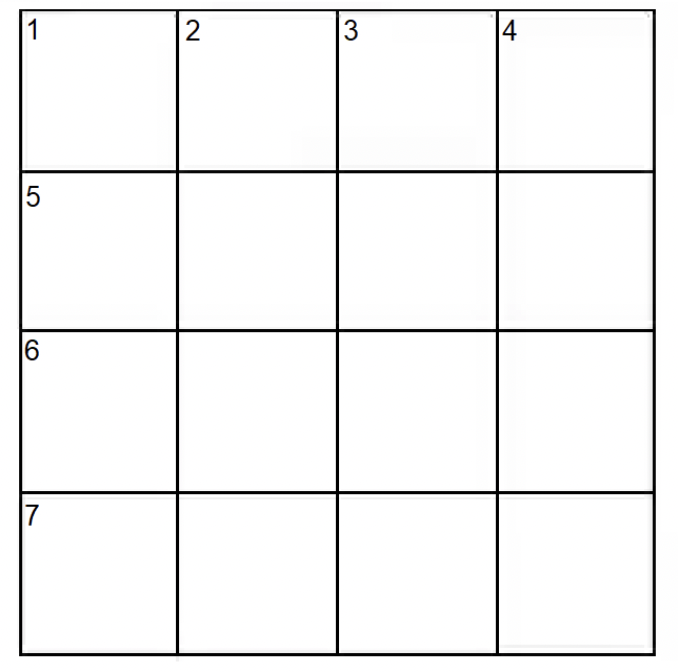}
    \captionof{figure}{$n=4,$ 8 clues}
    
   \end{minipage}\hfill
   \begin{minipage}{0.48\textwidth}
     \centering
      \includegraphics[width=0.35\linewidth]{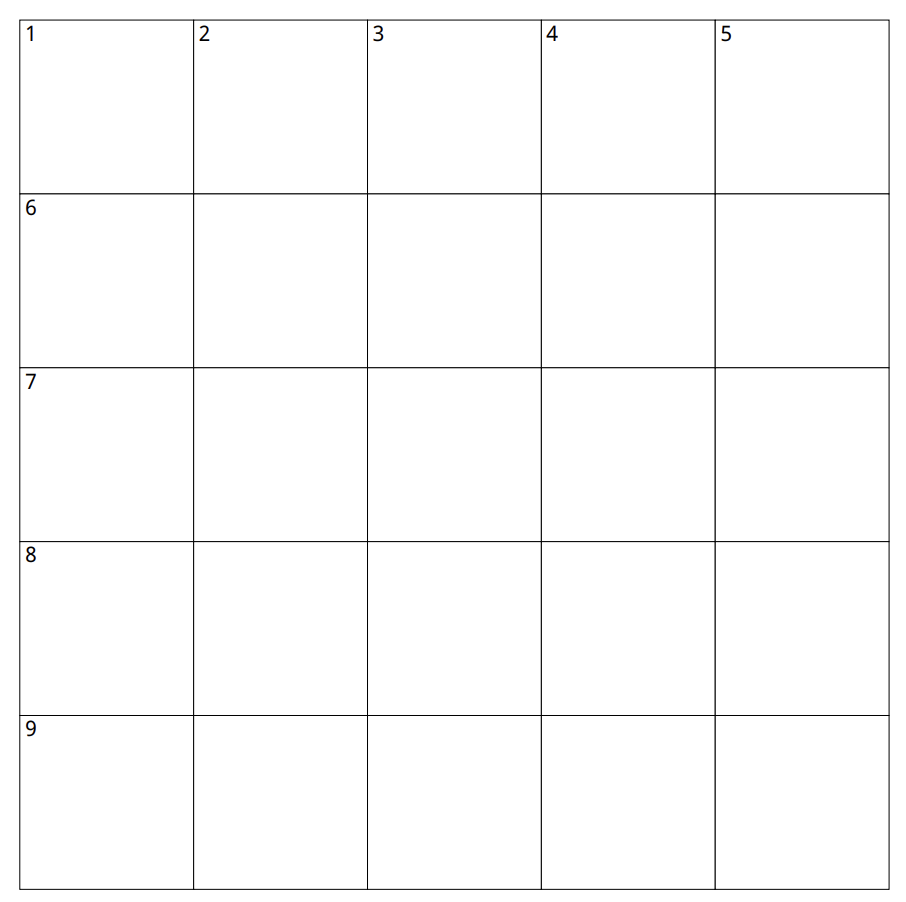}
       \captionof{figure}{$n=5,$ 10 clues}
    
   \end{minipage}
\end{figure}

\begin{figure}[!htb]
   \begin{minipage}{0.48\textwidth}
     \centering
      \includegraphics[width=0.35\linewidth]{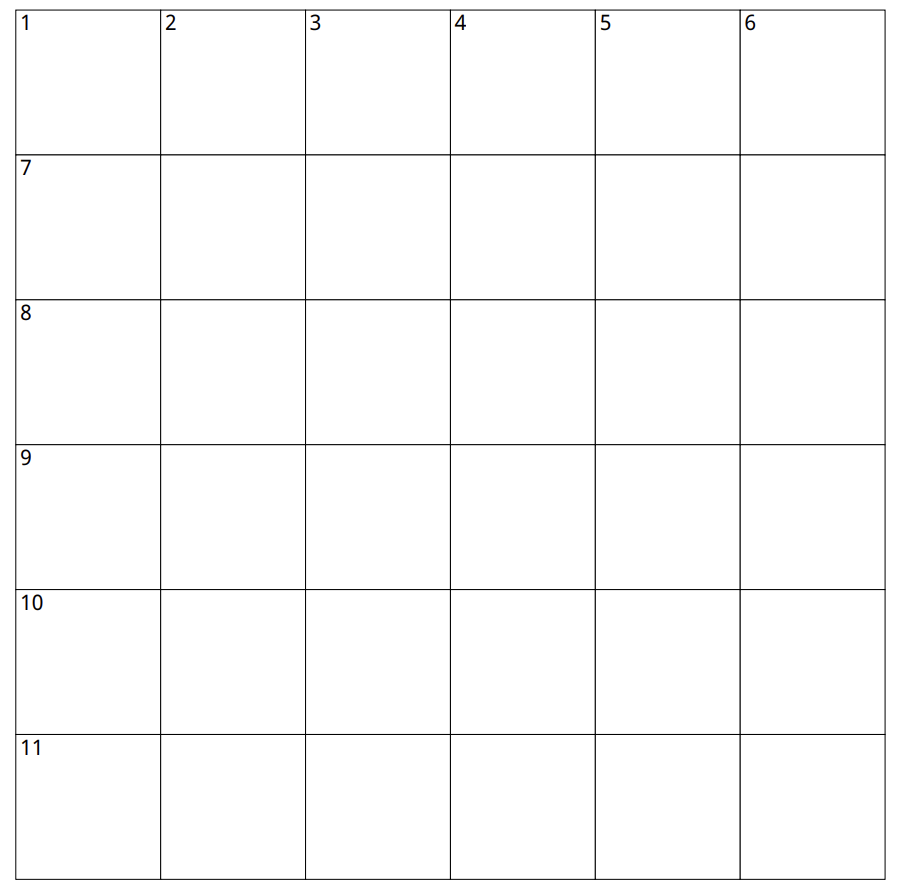}
    \caption{$n=6,$ 12 clues}
    
   \end{minipage}\hfill
   \begin{minipage}{0.48\textwidth}
     \centering
      \includegraphics[width=0.35\linewidth]{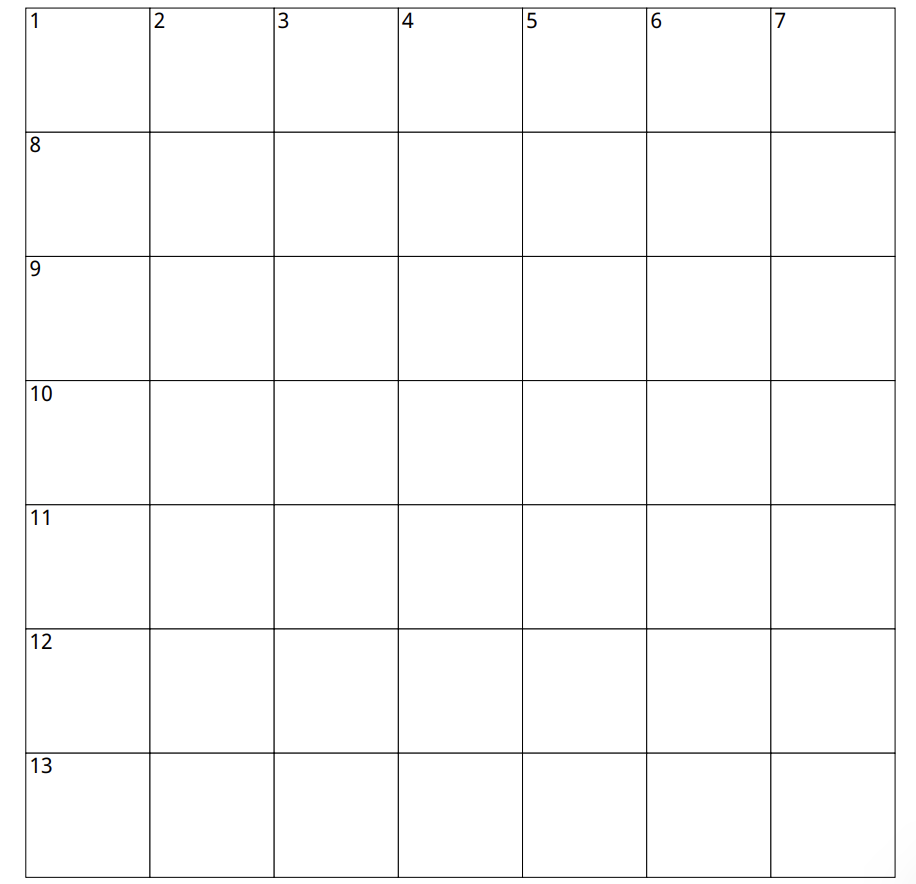}
    \caption{$n=7,$ 14 clues}
    
   \end{minipage}
\end{figure}

\begin{figure}[!htb]
   \begin{minipage}{0.48\textwidth}
     \centering
       \includegraphics[width=0.35\linewidth]{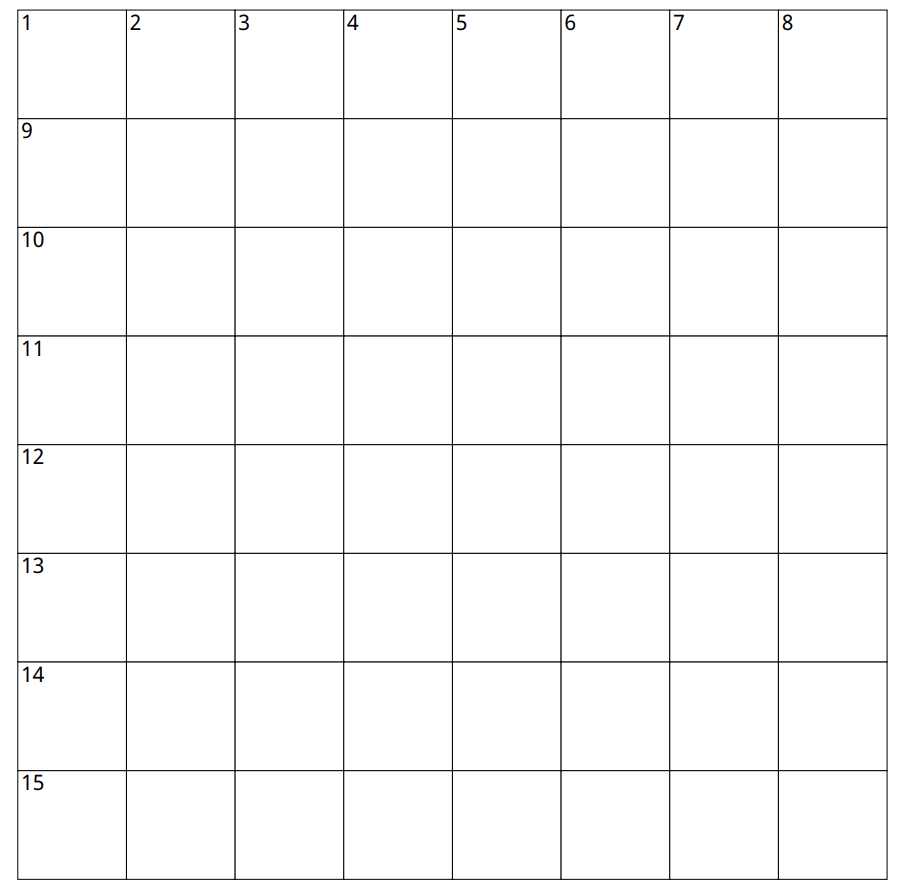}
    \caption{$n=8,$ 16 clues}
    
   \end{minipage}\hfill
   \begin{minipage}{0.48\textwidth}
     \centering
      \includegraphics[width=0.35\linewidth]{9.png}
    \caption{$n=9,$ 26 clues}
    \label{99}
   \end{minipage}
\end{figure}

\begin{figure}[!htb]
   \begin{minipage}{0.48\textwidth}
     \centering
       \includegraphics[width=0.35\linewidth]{10.png}
    \caption{$n=10,$ 36 clues}
    \label{1010}
   \end{minipage}\hfill
   \begin{minipage}{0.48\textwidth}
     \centering
      \includegraphics[width=0.35\linewidth]{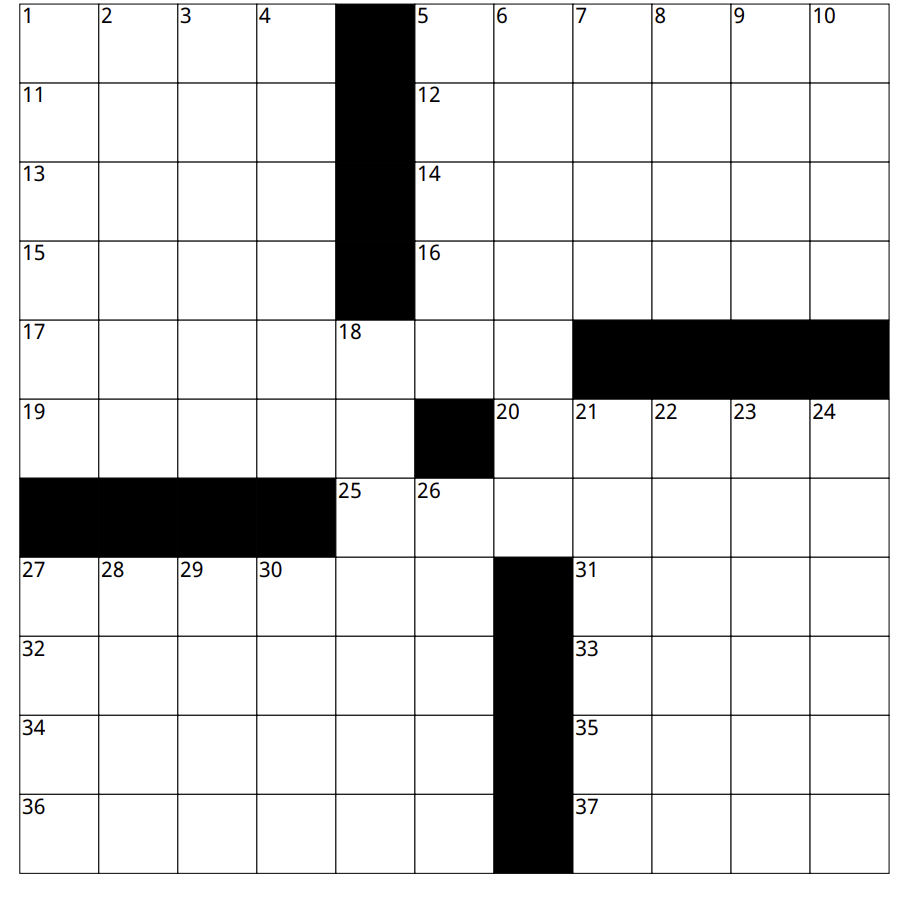}
    \caption{$n=11,$ 40 clues}
    
   \end{minipage}
\end{figure}

\begin{figure}[!htb]
   \begin{minipage}{0.48\textwidth}
     \centering
        \includegraphics[width=0.35\linewidth]{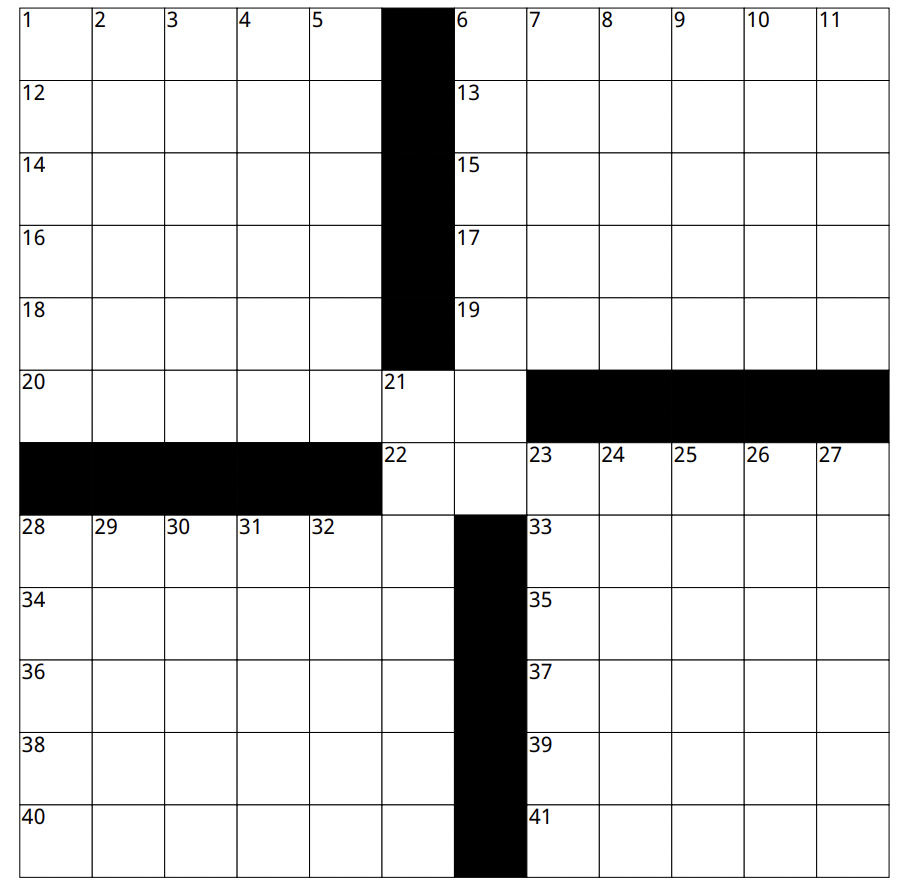}
    \caption{$n=12,$ 44 clues}
    
   \end{minipage}\hfill
   \begin{minipage}{0.48\textwidth}
     \centering
      \includegraphics[width=0.35\linewidth]{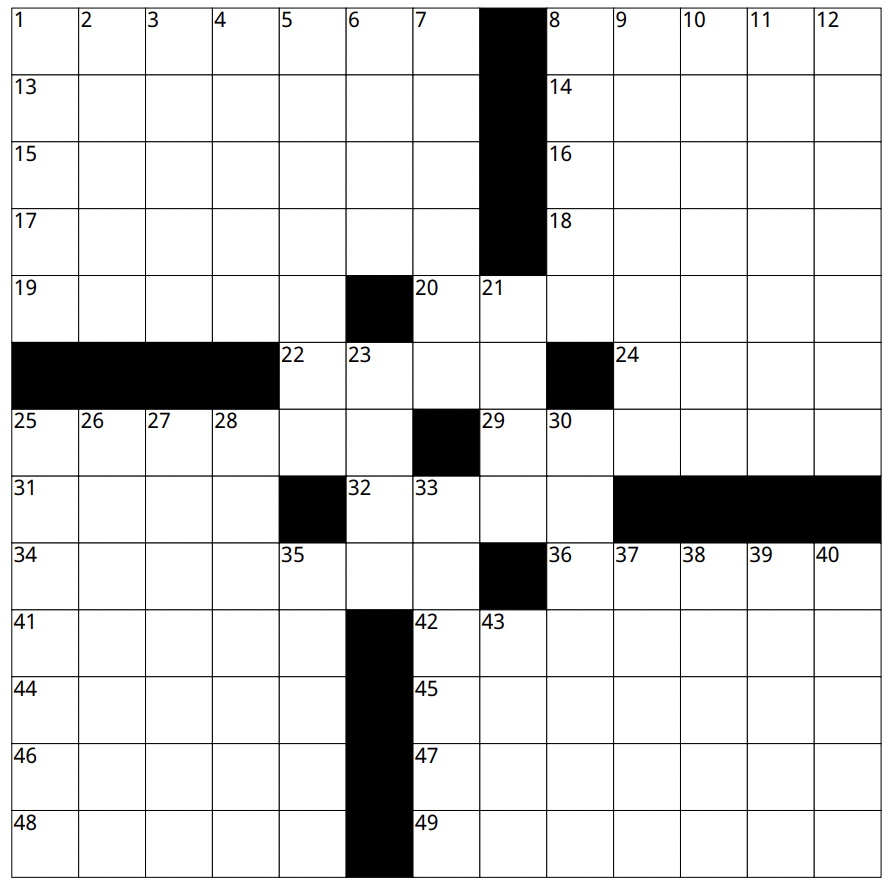}
    \caption{$n=13,$ 52 clues}
    
   \end{minipage}
\end{figure}

\begin{figure}[!htb]
   \begin{minipage}{0.48\textwidth}
     \centering
        \includegraphics[width=0.35\linewidth]{14.png}
    \caption{$n=14,$ 64 clues}
    \label{1414}
   \end{minipage}\hfill
   \begin{minipage}{0.48\textwidth}
     \centering
      \includegraphics[width=0.35\linewidth]{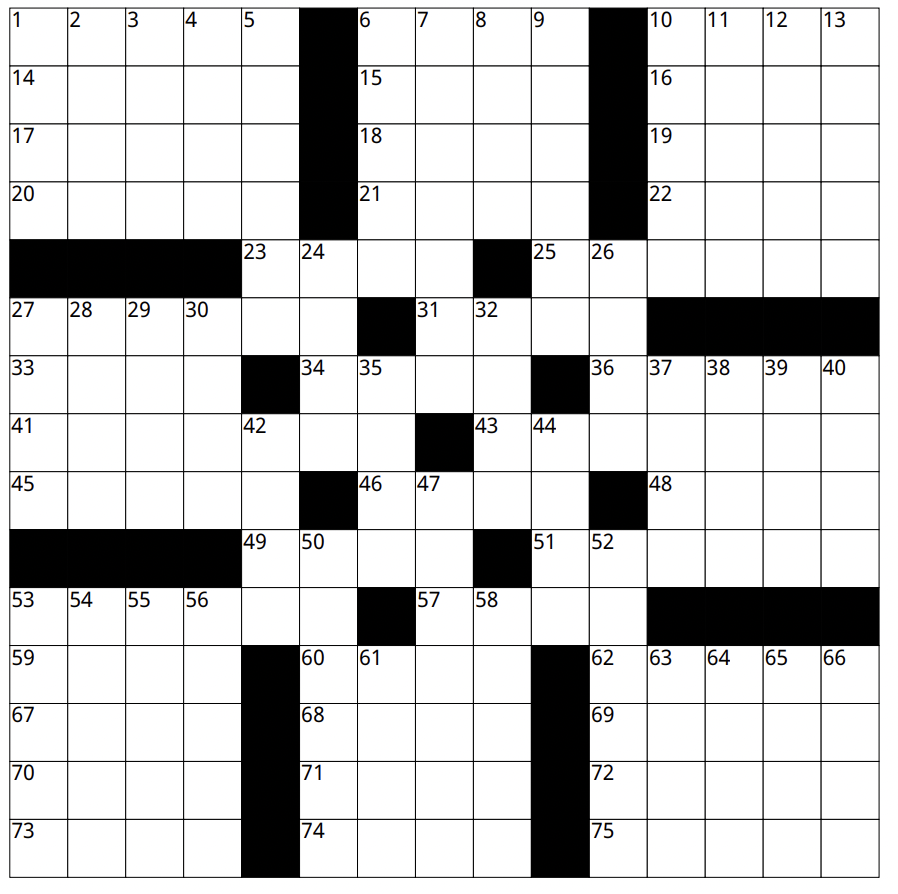}
    \caption{$n=15,$ 80 clues}
    
   \end{minipage}
\end{figure}

\begin{figure}[!htb]
   \begin{minipage}{0.48\textwidth}
     \centering
      \includegraphics[width=0.35\linewidth]{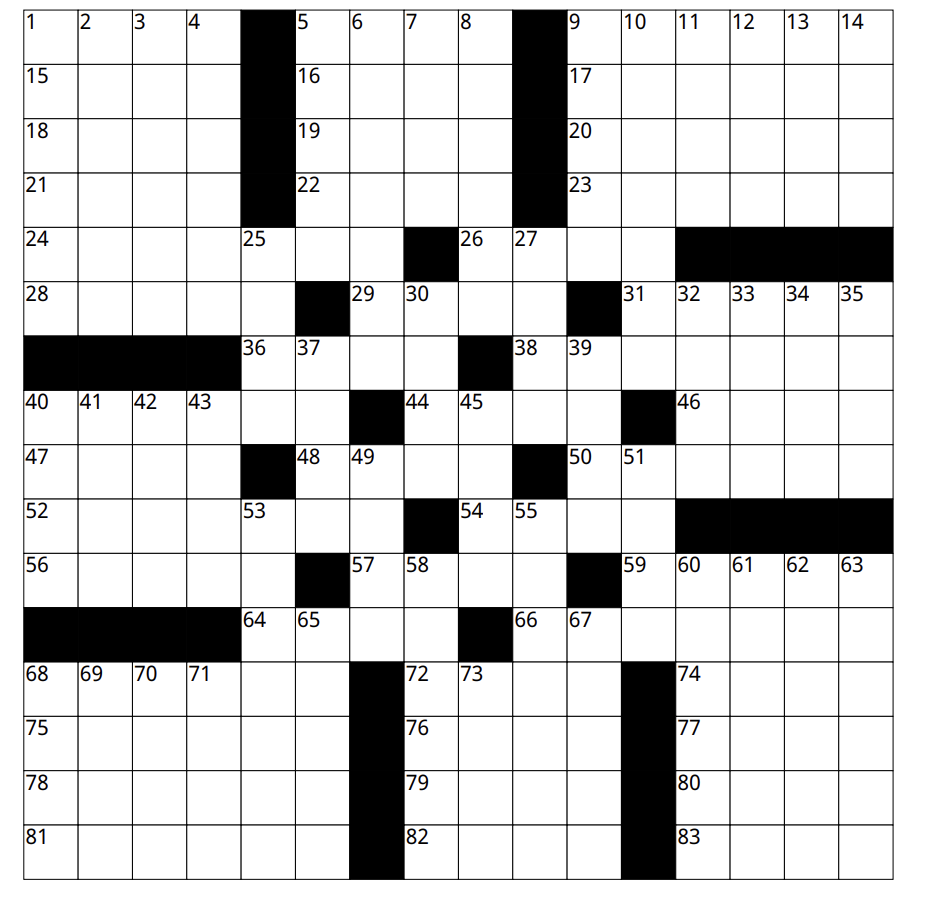} \par
    \caption{$n=16,$ 88 clues}
    
   \end{minipage}\hfill
   \begin{minipage}{0.48\textwidth}
     \centering
     \includegraphics[width=0.35\linewidth]{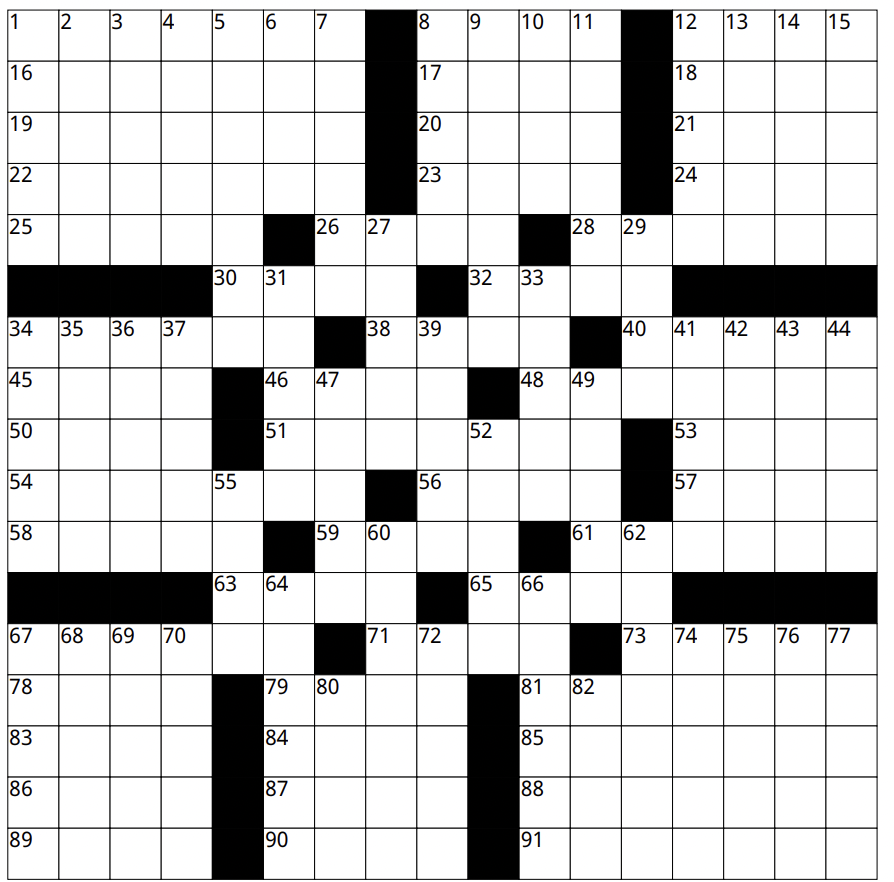} \par
    \caption{$n=17,$ 96 clues}
    
   \end{minipage}
\end{figure}

\begin{figure}[!htb]
   \begin{minipage}{0.48\textwidth}
     \centering
      \includegraphics[width=0.35\linewidth]{18.png} \par
    \caption{$n=18,$ 108 clues}
    
   \end{minipage}\hfill
   \begin{minipage}{0.48\textwidth}
     \centering
     \includegraphics[width=0.35\linewidth]{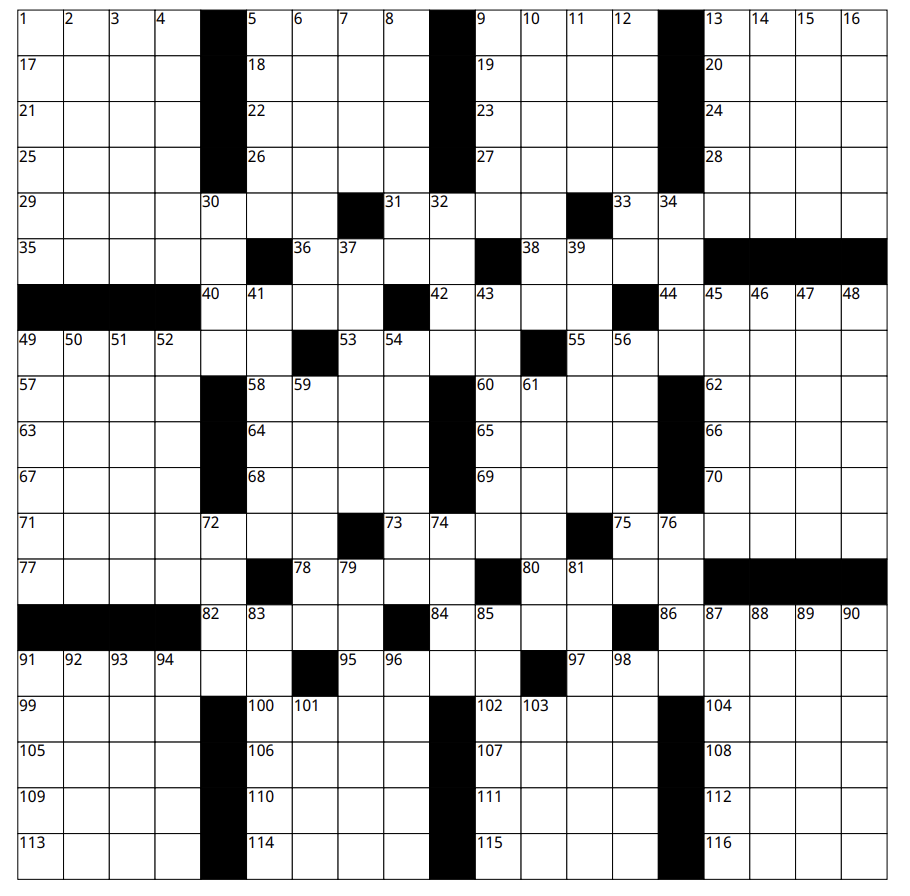} \par
    \caption{$n=19,$ 122 clues}
    
   \end{minipage}
\end{figure}

\begin{figure}[!htb]
   \begin{minipage}{0.48\textwidth}
     \centering
       \includegraphics[width=0.35\linewidth]{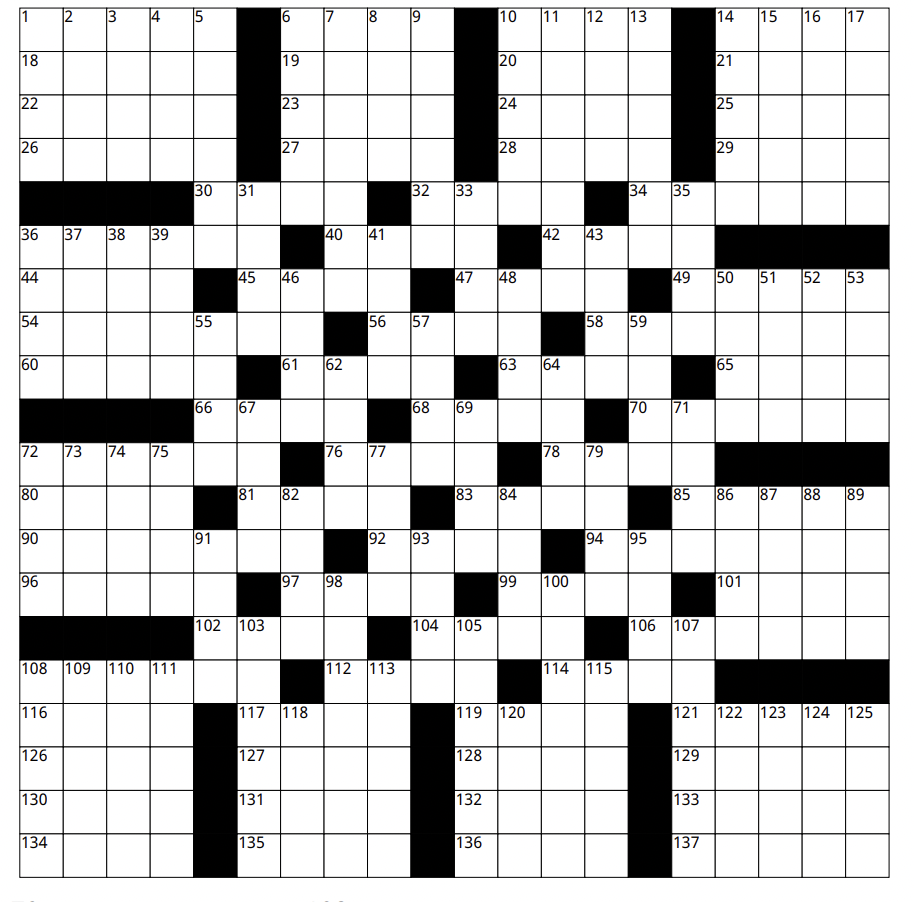} 
    \caption{$n=20,$ 144 clues}
    
   \end{minipage}\hfill
   \begin{minipage}{0.48\textwidth}
     \centering
     \includegraphics[width=0.35\linewidth]{21.png} 
    \caption{$n=21,$ 156 clues}
    
   \end{minipage}
\end{figure}

\begin{figure}[!htb]
   \begin{minipage}{0.48\textwidth}
     \centering
       \includegraphics[width=0.35\linewidth]{22.png} 
    \caption{$n=22,$ 166 clues}
    
   \end{minipage}\hfill
   \begin{minipage}{0.48\textwidth}
     \centering
     \includegraphics[width=0.35\linewidth]{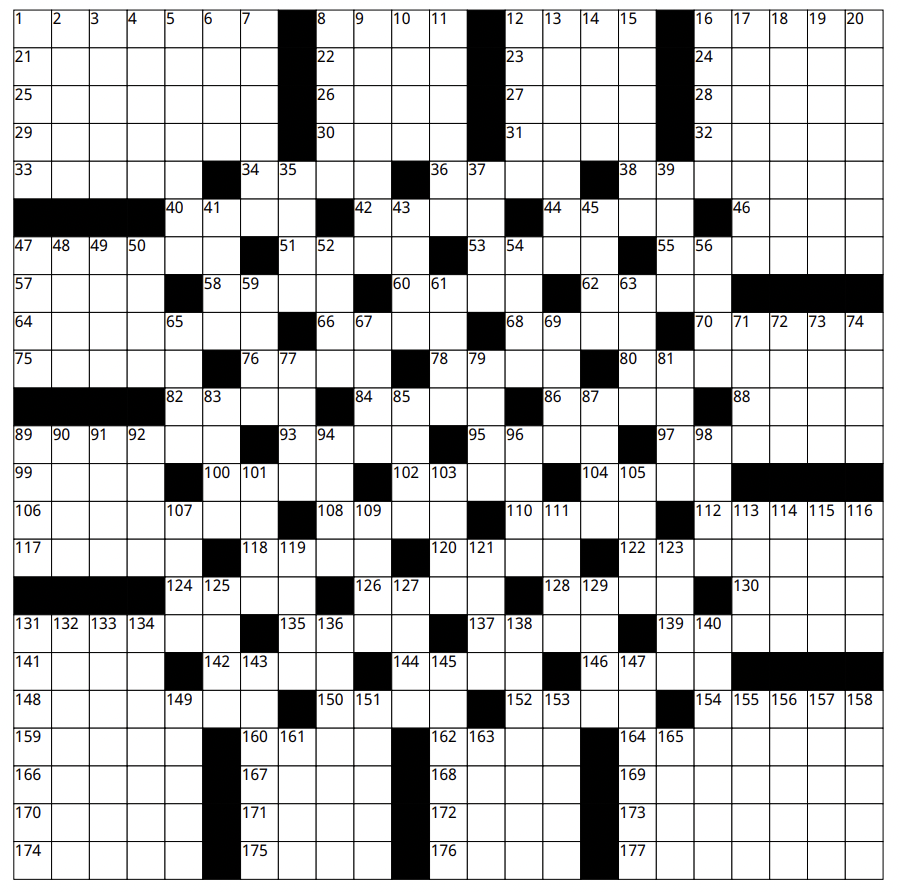} 
    \caption{$n=23,$ 184 clues}
    
   \end{minipage}
\end{figure}

\begin{figure}[!htb]
\begin{minipage}{0.48\textwidth}
 \centering
\includegraphics[width=0.35\linewidth]{24.png} 
    \caption{$n=24,$ 200 clues}
 
   \end{minipage}\hfill
   \begin{minipage}{0.48\textwidth}
     \centering
  \includegraphics[width=0.35\linewidth]{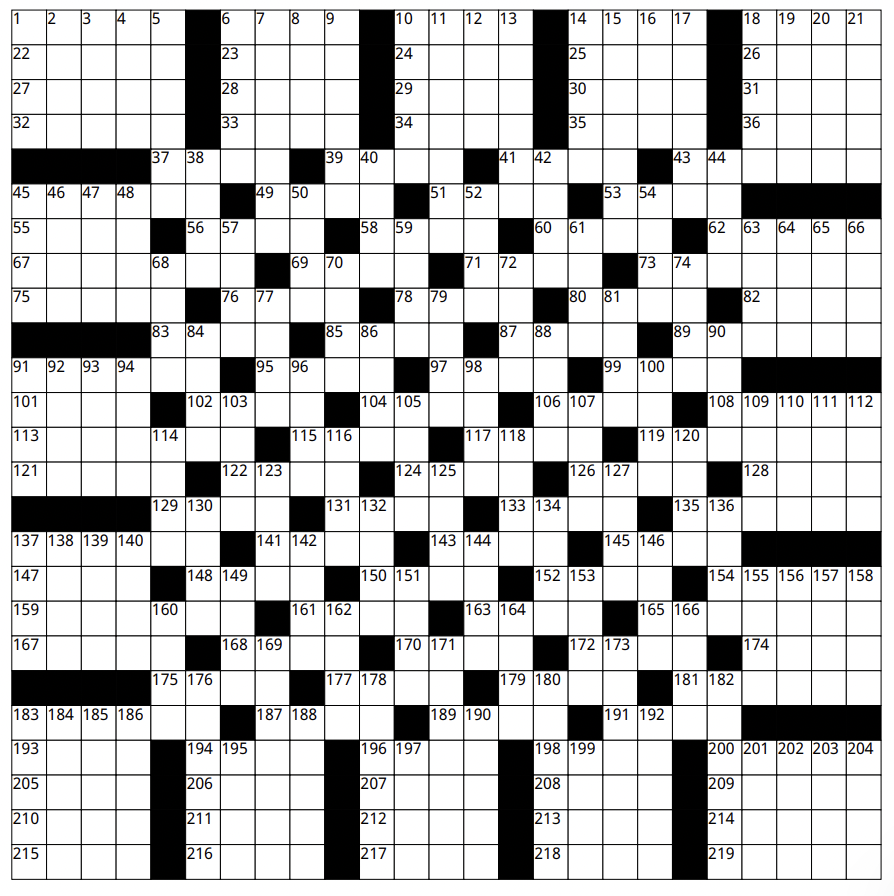} 
    \caption{$n=25,$ 228 clues}
    
   \end{minipage}
\end{figure}

\begin{figure}[!htb]
   \begin{minipage}{0.48\textwidth}
     \centering
       \includegraphics[width=0.35\linewidth]{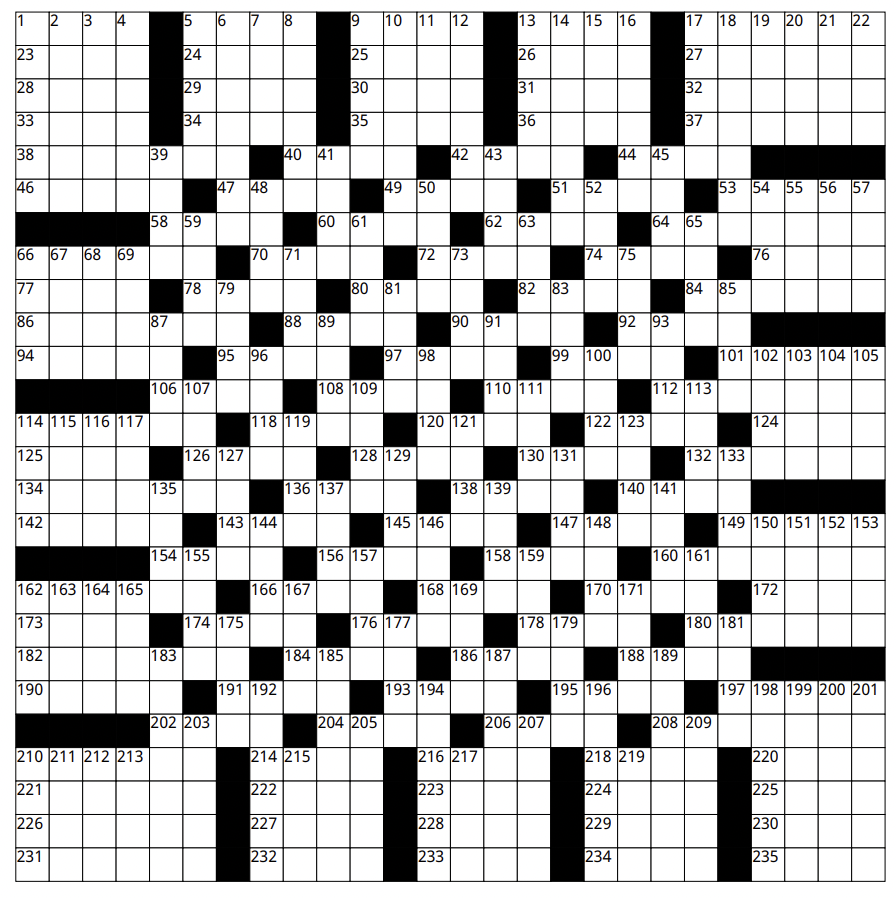} 
    \caption{$n=26,$ 244 clues}
    
   \end{minipage}\hfill
   \begin{minipage}{0.48\textwidth}
     \centering
     \includegraphics[width=0.35\linewidth]{27.png} 
    \caption{$n=27,$ 258 clues}
    
   \end{minipage}
\end{figure}

\begin{figure}[!htb]
   \begin{minipage}{0.48\textwidth}
     \centering
       \includegraphics[width=0.35\linewidth]{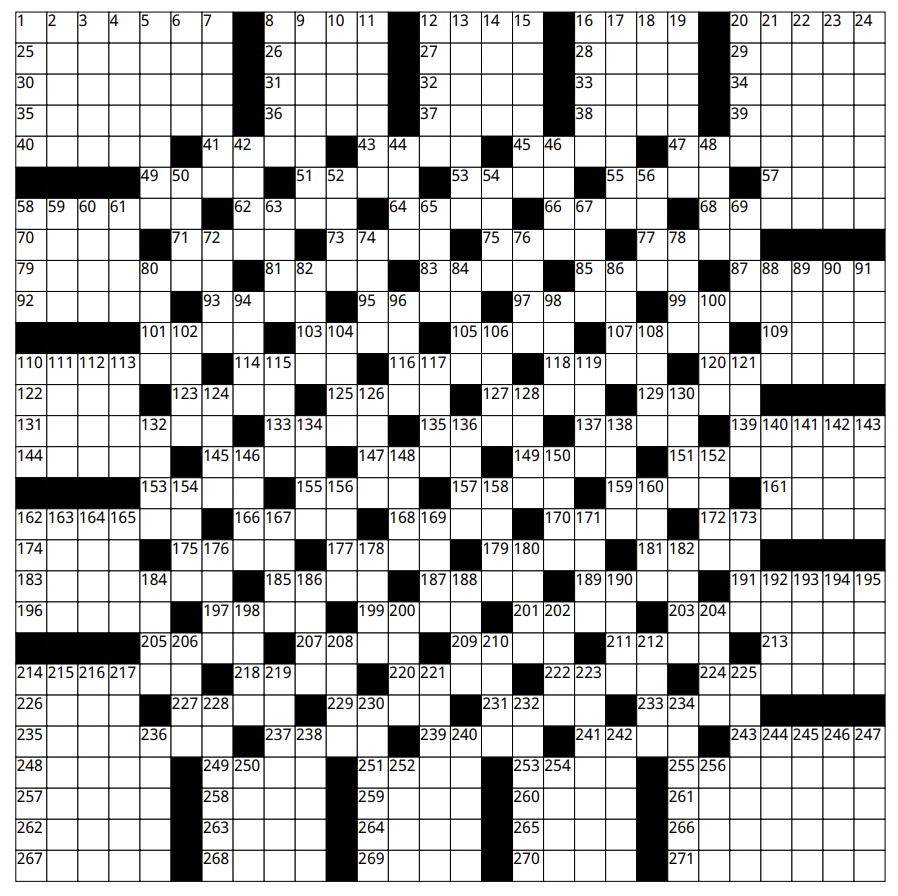} 
    \caption{$n=28,$ 280 clues}
    
   \end{minipage}\hfill
   \begin{minipage}{0.48\textwidth}
     \centering
      \includegraphics[width=0.35\linewidth]{29.png} 
    \caption{$n=29,$ 300 clues}
    
   \end{minipage}
\end{figure}

\begin{figure}[!htb]
   \begin{minipage}{0.48\textwidth}
     \centering
      \includegraphics[width=0.35\linewidth]{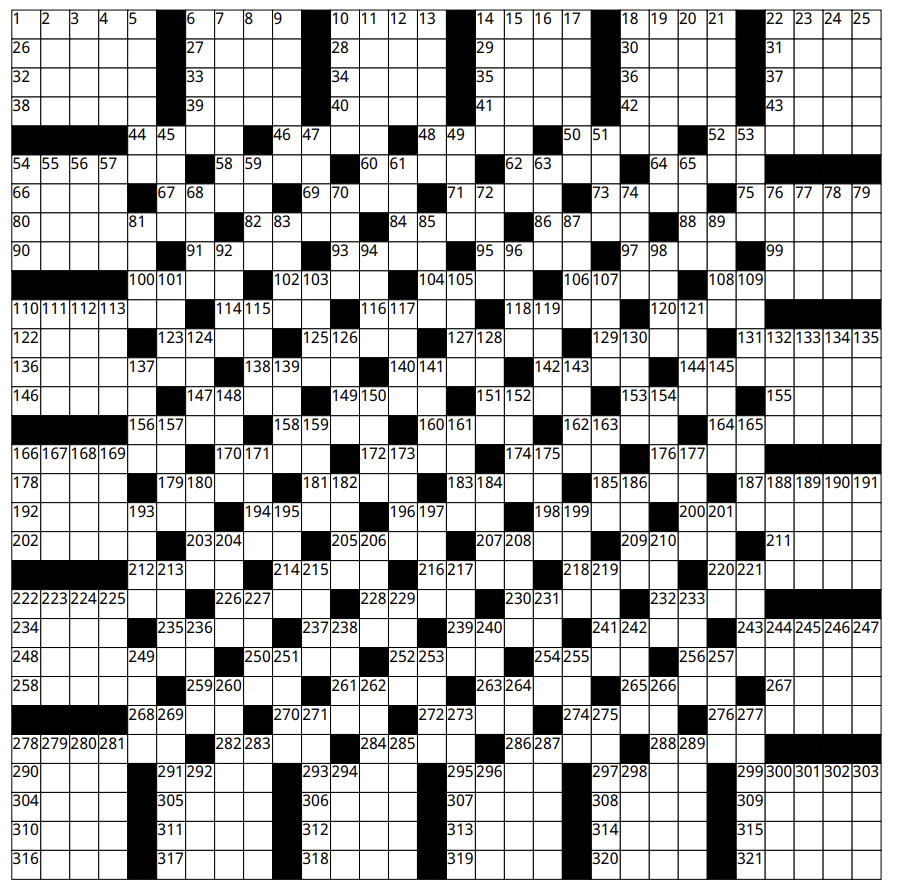} 
    \caption{$n=30,$ 332 clues}
    
   \end{minipage}\hfill
\end{figure}

\clearpage
\bibliographystyle{plain}
\bibliography{maximal_crosswords.bib}
\end{document}